\documentclass[leqno]{amsart}

\usepackage[foot]{amsaddr}
\usepackage{csquotes}

\usepackage{color}
\usepackage{geometry}
\usepackage{easyReview} 
\usepackage{subfiles}
\usepackage{graphicx}
\usepackage[section]{placeins} 
\usepackage{hyperref}
\usepackage{stmaryrd}
\usepackage{amsmath,amssymb,amsfonts,amsthm,textcomp}
\usepackage{mathtools}
\usepackage{tensor}
\usepackage{multicol}
\usepackage{subfig}
\usepackage{csquotes}
\usepackage{stackrel}
\usepackage{tikz}
\usepackage{tikz-3dplot}
\usetikzlibrary{3d}   
\usepackage{pgfplots}
\usepgfplotslibrary{fillbetween}
\pgfplotsset{compat=newest}
\usetikzlibrary{decorations.pathmorphing}
\usetikzlibrary{plotmarks}
\usetikzlibrary{arrows}
\usetikzlibrary{shapes.misc, positioning}
\usetikzlibrary{arrows,shapes,positioning}
\usetikzlibrary{decorations.markings}
\usetikzlibrary{arrows.meta,bending}
\usepackage{tikz}
\usepackage{pgfplots}
\pgfplotsset{compat=1.18}
\usetikzlibrary{arrows.meta}
\usepgfplotslibrary{fillbetween}

\pgfplotsset{
	myaxis/.style={
		thick,
		scale only axis,
		axis lines=middle,
		inner axis line style={-Triangle},
		axis on top,
		ticks=none,           
	},
}
\tikzset{
	obsfill/.style={fill=gray!25},
	obsline/.style={black, thick},
	hideline/.style={black, thick, dashed},
}
\usepackage{siunitx}
\usepackage{soul}
\usepackage{pdfpages}
\usepackage{dsfont}
\usepackage{bm}       
\usepackage{booktabs}
\usepackage{algorithm2e}
\usepackage{float}
\usepackage{enumerate}
\graphicspath{
              {Images/}
             }

\usepackage{mathtools}

\definecolor{ForestGreen}{rgb}{0.13, 0.55, 0.13}
\definecolor{BrickRed}{rgb}{0.80, 0.25, 0.33}
\definecolor{BlueViolet}{rgb}{0.54, 0.17, 0.89}
\definecolor{Plum}{rgb}{0.56, 0.27, 0.52}
\definecolor{Orange}{rgb}{1.0, 0.55, 0.0}

\tikzset{ForestGreen/.style={draw=ForestGreen,fill=ForestGreen}}
\tikzset{FGLabel/.style={text=ForestGreen}}
\tikzset{FGLine/.style={draw=ForestGreen,thick}}

\newfont{\tenbfsl}{cmbxti9 scaled 1200}
\newfont{\tenbbb}{msbm10}
\newfont{\svnbbb}{msbm8}

\DeclareMathOperator*{\argmin}{arg\,min}

\newcommand{\lrang}[2]{\left\langle {#1},\, {#2} \right\rangle}

\newtheorem{lem}{Lemma}

\newtheorem{theo}{Theorem}[section]      
\newtheorem{propo}[theo]{Proposition}

\newtheorem{defi}[theo]{Definition}
\newtheorem{exam}[theo]{Example}

\newtheorem{rem}[theo]{Remark}

\newtheorem*{exam*}{Example}             

\newcommand{\mcal}[1]{\mathcal{#1}}
\newcommand{\mb}[1]{\mathbb{#1}}
\newcommand{\mf}[1]{\mathbf{#1}}
\newcommand{\R}{\mathbb{R}}
\newcommand{\F}{\mathbb{F}}
\newcommand{\E}{\mathbb{E}}

\newcommand{\dl}{\mathrm{d}}
\newcommand{\ve}{\varepsilon}

\DeclareMathOperator{\diag}{diag}
\DeclareMathOperator{\trace}{tr}

\begin{document}

\title[PMP on the Belief Space for Continuous-Time Optimal Control with Discrete Observations]{A Pontryagin Maximum Principle on the Belief Space for Continuous-Time Optimal Control with Discrete Observations}
\author[C. Bayer, S. Ben Naamia, E. von Schwerin, R. Tempone]{Christian Bayer$^{1}$, Saifeddine ben Naamia$^{\#,2}$,  Erik von Schwerin$^{3}$ \& Ra\'{u}l Tempone$^{2,3,4}$}
\address{$^{1}$Weierstrass Institut (WIAS), Berlin, Germany,}
\address{$^2$Department of Mathematics, RWTH Aachen University, Geb\"{a}ude-1953 1.OG, Pontdriesch 14-16, 161, 52062 Aachen, Germany}
\address{$^3$King Abdullah University of Science \& Technology (KAUST), Computer, Electrical and Mathematical Sciences \& Engineering Division (CEMSE), Thuwal 23955-6900, Saudi Arabia}
\address{$^4$Alexander von Humboldt Professor in Mathematics for Uncertainty Quantification, RWTH Aachen University, Germany}
\email{$^\#$Bennaamia@ssd.rwth-aachen.de}


\date{\today}

\begin{abstract}
\noindent

We study a continuous time stochastic optimal control problem under partial observations that are available only at discrete time instants. This hybrid setting, with continuous dynamics and intermittent noisy measurements, arises in applications ranging from robotic exploration and target tracking to epidemic control. We formulate the problem on the space of beliefs (information states), treating the controller's posterior distribution of the state as the state variable for decision making. On this belief space we derive a Pontryagin maximum principle that provides necessary conditions for optimality. The analysis carefully tracks both the continuous evolution of the state between observation times and the Bayesian jump updates of the belief at observation instants.

A key insight is a relationship between the adjoint process in our maximum principle and the gradient of the value functional on the belief space, which links the optimality conditions to the dynamic programming approach on the space of probability measures. The resulting optimality system has a prediction and update structure that is closely related to the unnormalised Zakai equation and the normalised Kushner-Stratonovich equation in nonlinear filtering.

Building on this analysis, we design a particle based numerical scheme to approximate the coupled forward (filter) and backward (adjoint) system. The scheme uses particle filtering to represent the evolving belief and regression techniques to approximate the adjoint, which yields a practical algorithm for computing near optimal controls under partial information. The effectiveness of the approach is illustrated on both linear and nonlinear examples and highlights in particular the benefits of actively controlling the observation process.

\end{abstract}

\maketitle
\section*{Acknowledgments}
C. Bayer acknowledges support from DFG CRC/TRR 388 \enquote{Rough Analysis, Stochastic Dynamics and Related Fields} Project B03 and DFG individual grant number 497300407 \enquote{Recursive and sparse approximation in reinforcement learning with applications}.
S. Ben naamia acknowledges support from  the Deutsche Forschungsgemeinschaft (DFG, German Research Foundation) -- 333849990/GRK2379 (IRTG Hierarchical and Hybrid Approaches in Modern Inverse Problems).
This work was supported by the King Abdullah University of Science and Technology (KAUST) Office of Sponsored Research (OSR) under Award  No.~OSR-2019-CRG8-4033 and the Alexander von Humboldt Foundation.

\tableofcontents                        


%
\section{Introduction}
\label{sec:intro}

Many control systems operate under partial information: the controller cannot directly observe the full state of the system and instead has access only to partial, noisy measurements. Classical examples include navigation and tracking with intermittent sensor readings, robotic exploration with limited feedback from the environment, and epidemic control where infection states must be inferred from sparse testing data. In such settings, control actions can influence not only the state evolution but also the quality and timing of future information, a phenomenon known as the dual effect of control. The controller must balance regulating the system and probing it to gather information, and optimal strategies under partial observations are typically much more complex than in fully observed control problems.

Even in the linear quadratic Gaussian (LQG) setting, where a separation principle holds, the output feedback problem is already nontrivial: the optimal controller can be implemented by combining a Kalman filter with a full state LQR law. Outside this narrow regime, separation and certainty equivalence may fail. Witsenhausen's counterexample \cite{doi:10.1137/0306011}shows that an optimal controller for a seemingly simple problem with a single nonlinearity and partial information can differ drastically from any design based on naive separation ideas. In general, partially observed stochastic control, also known as optimal control with incomplete information, rarely admits closed form solutions and remains an active area of research.

\subsection{Literature review} \label{subsec:lit}

A classical way to handle partial observations is to reformulate the problem as a fully observed control problem on the space of beliefs (information states), that is, the conditional distribution of the state given the observation history. This idea goes back to early work on stochastic control and filtering and leads to a dynamic programming equation (Hamilton-Jacobi-Bellman, HJB) on an infinite dimensional space of probability measures; see for example \cite{Astrom1970,KumarVaraiya1986,BertsekasShreve1978,Bensoussan1992,FS93,BCD97}. While this formulation is conceptually clean and underlies the separation principle in the LQG case, solving the resulting measure valued HJB equation is computationally intractable except in special structures, such as finite state models or problems with finite dimensional sufficient statistics.

Another approach is based on stochastic maximum principles. Here one augments the state with the filtering dynamics and derives first order necessary conditions in terms of a forward state (or filter) and a backward adjoint process. Early contributions already imposed the filter dynamics explicitly in order to apply Pontryagin's principle under partial information \cite{BarasElliottKohlmann1989,Bensoussan1992}. Later, Tang \cite{Tang1998} established a general maximum principle for partially observed stochastic differential systems with correlated noise between state and observation. Very recently, Opper and Reich \cite{opper2025digitaltwinsmckeanpontryagincontrol} have developed a McKean–Pontryagin minimum principle for partially observed systems, combining ensemble Kalman filtering with a mean-field optimal control formulation for real-time digital twin applications. 
These works show that Pontryagin type conditions can be extended to partial information at the price of handling an infinite dimensional, measure valued state.

In recent years optimal control problems with distribution dependent dynamics have been analysed directly on spaces of probability measures endowed with the Wasserstein metric. In particular, Bonnet and Rossi \cite{bonnet2017} and Bonnet \cite{bonnet2019} develop second order maximum principles and weak optimal control for such systems, and provide tools to differentiate cost functionals on the space of probability measures. This point of view supplies a natural mathematical framework for our belief space formulation and justifies the use of variational derivatives with respect to the law of the state.

The present work is motivated by continuous time systems in which observations are available only at discrete instants. This leads to a hybrid structure: the physical state evolves continuously, whereas the information available to the controller is updated in jumps at observation times. Classical nonlinear filtering theory, including the Zakai equation and the Kushner-Stratonovich equation, covers continuous time observation streams; in contrast, the discrete time observation case requires a different treatment. A recent preprint by Bayer et al.\ \cite{bayer2024continuoustimestochasticoptimal} analyses such partially observed control problems via dynamic programming on the belief space and derives an HJB equation on an infinite dimensional space of probability measures. Our contribution is complementary: we develop a Pontryagin maximum principle on the same belief space and identify an adjoint process that can be related to the gradient of the value functional.

Another line of research seeks to mitigate the infinite dimensional nature of the belief by introducing finite dimensional memory or compressed statistics. Tottori and Kobayashi \cite{tottori2022memorylimited,tottori_forward-backward_2023} consider memory limited partially observable control, where the controller optimises with respect to a fixed finite dimensional information state updated by a prescribed compression map instead of the full posterior. They derive forward and backward systems that couple a Fokker-Planck equation for the state distribution with an HJB equation for the approximate belief. Our formulation is related in spirit: we introduce a finite dimensional memory state in the numerical scheme, but our theoretical optimality system is derived directly on the belief space with exact Bayesian updates.

Finally, there is a large literature on the separation principle, the dual effect, and active sensing under partial information. In linear Gaussian models with exogenous observations, the separation theorem \cite{Kalman1960,KalmanBucy1961,Wonham1968,Stengel1994} shows that optimal output feedback control can be implemented by combining a Kalman filter with a full state LQR controller. Outside this setting, certainty equivalence can fail, as illustrated by Witsenhausen's counterexample \cite{Witsenhausen1968}. When the observation law itself depends on the control policy, the dual effect destroys separation \cite{BarShalomTse1974,RameshSandbergBaoJohansson2011}, and the controller must trade off immediate control cost against information acquisition. This has motivated a growing body of work on active information gathering and controlled sensing in robotics and related areas \cite{cognetti2018,lauri2014,olivos2024}. Our work contributes to this broader theme by providing a belief space Pontryagin framework and a particle based algorithm that explicitly capture the coupling between control, observation design, and belief dynamics in continuous time with discrete observations.
\subsection{Our contribution}

Our contributions are fourfold.

First, we formulate a continuous time stochastic control problem with partial observations at discrete times in a general setting. The problem allows two kinds of control actions: (i) a continuous control $\alpha_t$ that affects the state dynamics between observations, and (ii) a discrete control $\beta_{t_n}$ applied at observation times that influences the observation process (for example by choosing sensor parameters or triggering measurements). To manage the information structure we introduce a finite dimensional memory state $Z_{t_n}$ updated at each observation by a prescribed compression mapping $\varphi_n$. This memory state summarises the past observations in finite dimension and makes it possible to describe controlled sensing and filtering in a single optimal control framework.

Second, we derive a Pontryagin maximum principle on the belief space. By lifting the problem to the space of probability measures we obtain necessary conditions in terms of a forward evolution for the belief and a backward adjoint process. Between observation times the belief satisfies a controlled Fokker-Planck equation and the adjoint satisfies a backward Kolmogorov type equation; at each observation instant both processes undergo a Bayesian jump. We give explicit jump conditions for the adjoint, including the extra term that accounts for the normalisation of the posterior. This structure mirrors the relationship between unnormalised and normalised nonlinear filtering equations and clarifies how control and information updates interact in continuous time with discrete observations.

Third, we establish a link between the adjoint variables and the value function in the belief space formulation. Under suitable convexity and differentiability conditions we show that the adjoint (costate) process can be identified with the gradient of the dynamic programming value functional on the space of probability measures, evaluated along the optimal trajectory. In particular, if $V_t(\mu)$ denotes the optimal value starting at time $t$ with current belief $\mu$ and $\bar\mu_t$ is the optimal belief trajectory, then the optimal adjoint $\bar U_t$ satisfies $\langle \bar U_t,\bar\mu_t\rangle = V_t(\bar\mu_t)$ and serves as an affine majorant of $V_t$ around $\bar\mu_t$. This relationship is the infinite dimensional counterpart of the classical envelope condition $p(t)=\nabla_x V(t,x)$ in fully observed optimal control and connects the Pontryagin maximum principle with the HJB equation in the partially observed setting.

Fourth, we design a numerical algorithm to compute locally optimal policies under partial observations and demonstrate its performance on linear and nonlinear examples. The method iteratively solves the forward and backward parts of the optimality system by simulation and function approximation. A particle filter is used to simulate the state and observation dynamics, producing an empirical approximation of the belief, while a parametric ansatz $U_t^{\theta}(x,\mf y^{[n]})$ (for example a polynomial or a neural network) is used for the adjoint. Given a current parameter $\theta$, the algorithm extracts policies $(\alpha_t,\beta_{t_n})$ by minimising belief weighted Hamiltonians, simulates many trajectories to estimate pathwise costs, and then regresses the parametric representation onto these costs to update $\theta$. To keep the regression problem tractable as the observation history grows, we condition the policy and value only on a sliding window of the most recent observations. The resulting scheme combines particle filtering, ideas from dynamic programming, and regression, and it can be applied to high dimensional and nonlinear partially observed problems.

\subsection{Organization of the paper}

Section~2 introduces the partially observed control model in continuous time with discrete observation updates. It derives the belief dynamics and presents the Pontryagin optimality system, including the jump conditions at observation times and the relation between the belief space value function and the adjoint. Section~3 describes the numerical approach: the particle representation of beliefs, the parametric ansatz for the adjoint or value function, and the iterative algorithm for policy optimisation, together with numerical experiments. Section~4 contains concluding remarks and discusses limitations and possible extensions. Technical proofs are collected in the appendix.

\subsection{Notation}

We briefly summarise the main notation and conventions used in the paper. The state of the system is denoted by $X_t\in\mathbb{R}^{d_x}$ at time $t\in[0,T]$ and evolves according to a controlled stochastic differential equation. The control has two components: a continuous control $\alpha_t\in\mathbb{R}^{d_\alpha}$ applied $\forall t \in [0,T]$  and a discrete control $\beta_{t_n}\in\mathbb{R}^{d_\beta}$ applied at observation times $t_n$. Observations $Y_{t_n}\in\mathbb{R}^{d_y}$ are received at discrete times $0<t_1<\dots<t_{N_o}<T$; no observations occur between these instants. The observation model at time $t_n$ has the form
\[
Y_{t_n}=h_n(t_n,X_{t_n},\beta_{t_n},\xi_n),
\]
where $\xi_n$ is the measurement noise.
For a fixed observation history we introduce the concatenation of all the available observations up to time $t \in [t_n,t_{n+1})$ denoted by $Y^{[n]}:=(Y_{t_1},\dots,Y_{t_n})$ and we write : 
\[
\mathcal{F}^Y_t=\sigma\bigl(Y^{[\nu(t)]}\bigr)
\]
for the observation filtration, where $\nu(t)$ is the index of the last observation received by time $t$. Admissible controls are adapted to the appropriate filtration, as specified in Section~2. The conditional distribution of $X_t$ given $\mathcal{F}^Y_t$ is called the \emph{belief} and is denoted by
\[
\mu_t(\cdot)=\mathbb{P}(X_t\in\cdot\mid\mathcal{F}^Y_t).
\]
We write $\mathcal{P}(\mathbb{R}^{d_x})$ for the space of Borel probability measures on $\mathbb{R}^{d_x}$. For a measurable function $U$ and a measure $\mu$ we use the pairing
\[
\langle U,\mu\rangle := \int_{\mathbb{R}^{d_x}} U(x)\,\mu(dx).
\]
Expectation with respect to all sources of randomness is denoted by $\mathbb{E}[\cdot]$. Other notation is introduced locally and used consistently with these conventions.


\section{Partially Observable Stochastic Optimal Control with Discrete Observations}
\label{sec:POSOC}

\subsection{Partially observed stochastic optimal control (POSOC) formulation}

In the partially observed setting, we only access partial, noisy measurements of the state process $X_t$.  
Observations occur at discrete times $0=t_0<t_1<\dots<t_{N_o}<T=:t_{N_o+1}$ and are modeled by the process $\{Y_{t_n}\}_{n=1}^{N_{o}}$.  
We consider the controlled dynamics
\begin{equation}
	\label{eq:POSOC_dynamics}
	\begin{aligned}
		\dl X_t &= b_t\!\left(X_t,\alpha_t\right)\,\dl t 
		+ \sigma_t\!\left(X_t,\alpha_t\right)\,\dl W_t,&& \forall t\in[0,T]\\
		Y_{t_n} &= h_n\!\left(t_n,X_{t_n},\beta_{t_n},\xi_n\right), && n=1,\dots,N_o,\\
		Z_{t_n} &= \phi_n\left(Z_{t_{n-1}},Y_{t_n}\right) \\
		X_0 &\sim \mu_0\in\mcal P(\mb{R}^{d_x}),
	\end{aligned}
\end{equation}
where $b:[0,T]\times\mb{R}^{d_x}\times\mb{R}^{d_\alpha}\to\mb{R}^{d_x}$,  
$\sigma:[0,T]\times\mb{R}^{d_x}\times\mb{R}^{d_\alpha}\to\mb{R}^{d_x\times d_w}$, and  
$W$ is a $d_w$–dimensional standard Wiener process.  
The observation $Y_{t_n}\in\mb{R}^{d_y}$ is generated by an observation function $h_n$ and the independent discrete driving noise for the measurement processes $\{\xi_n\}$. Let $u_t=(\alpha_t,\beta_t)$ denote the control functions.
We augment the model with a finite-dimensional memory state $Z_{t_n}\in\mathbb{R}^{d_z}$ that is updated only at observation times by measurable maps
$$
Z_{t_n} \;=\; \phi_n\!\left(Z_{t_{n-1}},\, Y_{t_n}\right), 
\qquad n=1,\dots,N_o, 
\quad Z_{t_0}=z_0\in\mathbb{R}^{d_z}.
$$
We interpret $\phi_n$ as a \emph{compression/feature-extraction} rule that selectively retains information from the new observation $Y_{t_n}$ and the previous compressed memory $Z_{t_{n-1}}$. In general, $\phi_n$ may itself be chosen by the designer (and thus viewed as a controllable component of the sensing/memory architecture). In this work, we fix $\{\phi_n\}_{n=1}^{N_o}$ in an ad hoc manner to illustrate the framework, while keeping the formulation broad to enable future \emph{co-design} of sensing, compression, and control. We also assume  knowledge of $b,\sigma_t ,h_n $ and $\mu_0$.
 Let $\mf Y^{[n]}:=\{Y_{t_1},\dots,Y_{t_n}\}$ and for $ t \in [0,T]$, set 
$$
\nu(t) := \max\{k \in \{1,\dots,N_o\} : t_k \le t\},
$$
with the convention \(\max\emptyset=0\) and \(\mf{Y}^{[0]}:=\emptyset\) and define the filtrations by the  \emph{$\sigma$-algebras} 
\begin{align}
	\mb F_t^Y &:= \sigma\!\big(\mf Y^{[\nu(t)]}\big),&& t\in[0,T],
	\label{eq:filtration-Y}\\
	\mb F_t^X&:=\sigma(X_s:0\le s\le t), && t\in[0,T],
	\label{eq:filtration-X}\\
	\mb F_t^{X,Y} &:= \mb F_t^X \vee \mb F_t^Y = \sigma\big(\mathcal{F}^X_t \cup \mathcal{F}^Y_t\big),&& t\in[0,T].
	\label{eq:joint-filtration}
\end{align}
so that $\mb F_t^Y$ is piecewise constant on $[0,T]$ with jumps at $t_n, n=1,\dots,N_o$.

\subsubsection{Controls and admissibility}
\emph{Throughout, we adopt a closed-loop (feedback) formulation: admissible controls are progressively measurable with respect to the relevant information filtration (e.g., ${\mb F^Y_t}$ or ${\mb F^X_t}$) and are implemented as functionals of the available observations/state, rather than as open-loop time plans.}

We consider continuous controls $\alpha_t\in\mb{R}^{d_\alpha}$ and discrete controls $\beta_{t_n}\in\mb{R}^{d_\beta}$:
\begin{defi}[Admissible controls]
	\label{def:admissible}
	Let $\mcal U[0,T]$ be the set of all $u=(\alpha,\beta)$ such that
	\begin{itemize}
		\item $\alpha$ is $\{\mb F_t\}$–adapted on $[0,T]$;
		\item $\beta$ is piecewise constant, left–continuous with jumps only at $\{t_n\}$, and $\beta_{t_n}$ is $\mb F_{t_n^-}$–measurable (hence predictable);
		\item $u_t=(\alpha_t,\beta_t)\in \mb{R}^{d_\alpha}\times\mb{R}^{d_\beta}$ for all $t$;
	\end{itemize}
	where  $\mb F_t$ can be $\mb F_t^Y$ or $\mb F_t^X$.
\end{defi}

\begin{defi}[Three control classes]
	\label{rem:control-classes}
	For clarity we distinguish between three different  information structures that we can use to pick the control. For $t\in[0,T]$ define
	\[
	\begin{aligned}
		\mcal U^{\mathrm{PO}}[t,T]
		&:= \Big\{u=(\alpha,\beta)\in \left(\mcal A, \mcal B \right):\ \alpha \text{ is }\mb F^Y\text{-adapted on }[t,T],\ 
		\beta_{t_n}\text{ is }\mb F^Y_{t_n^-}\text{-measurable for }t_n\ge t\Big\},\\[0.25em]
		\mcal U^{\mathrm{FO}}[t,T]
		&:= \Big\{u=(\alpha,\beta)\in  \left(\mcal A, \mcal B \right):\ \alpha \text{ is }\mb F^X\text{-adapted on }[t,T],\ 
		\beta_{t_n}\text{ is }\mb F^X_{t_n^-}\text{-measurable for }t_n\ge t\Big\},\\[0.25em]
		\mcal U^{\mathrm{R}}[t,T]
		&:= \Big\{u=(\alpha,\beta)\in  \left(\mcal A, \mcal B \right):\ \alpha \text{ is }\mb F^{X,Y}\text{-adapted on }[t,T],\ 
		\beta_{t_n}\text{ is }\mb F^{X,Y}_{t_n^-}\text{-measurable for }t_n\ge t\Big\},
	\end{aligned}
	\]
	where $ \left(\mcal A, \mcal B \right)$ are the set of  constraints that the controls $\left(\alpha_s,\beta_s\right)$ have to satisfy for  $ s \in [t,T]$. The \emph{partially observed} class $\mcal U^{\mathrm{PO}}$ is the feasible set for the POSOC problem, while the \emph{fully observed} class $\mcal U^{\mathrm{FO}}$ (state–feedback admissible) is used for the lower envelope. The auxiliary class $\mcal U^{\mathrm{R}}$ is found useful later in the proof of Proposition \ref{prop:envelope-ineq}.
\end{defi}

\subsubsection{Objective function}
Given running, impulse, and terminal costs
\[
f_t:\mb{R}^{d_x}\times\mb{R}^{d_\alpha}\to\mb{R},\quad
c_{t_n}:\mb{R}^{d_x}\times\mb{R}^{d_\beta}\to\mb{R},\quad
g:\mb{R}^{d_x}\to\mb{R},
\]
the expected cost of $u\in\mcal U[0,T]$ is
\begin{equation}
	\label{eq:POSOC_cost_functional-X}
	J(u)
	= \mb E\!\left[
	\int_0^T f_t\!\left(X_t,\alpha_t\right)\,\dl t
	+ \sum_{n=1}^{N_o} c_{t_n}\!\left(X_{t_n^-},\beta_{t_n}\right)
	+ g(X_T)
	\right],
\end{equation}
and the POSOC problem is $\min_{u\in\mcal U} J(u)$.

\begin{exam}[POSOC-LQG]
	In this LQG example, we assume access only to discrete, noisy measurements of the state process $X_t$. These observations are represented by the stochastic process  $Y_{t}$ , sampled at discrete time points $t_n, n=1,\dots,N_o$.
	\begin{equation}
		\label{eq:POSOC_LQG_observation}
		\begin{aligned}
			\dl X_t &= \left(A \, X_t + B \, \alpha_t \right) \, \dl t \,+\, \sigma \, \dl W_t,\\
			Y_{t_n} &= C \, X_{t_n} + \diag(\beta_{t_n}) \xi_n, \quad n=1,\dots,N_o,\\
			X_0 &\sim \mcal{N}(m_0,\Sigma_0),
		\end{aligned}
	\end{equation}
	where $\xi_n \sim \mcal{N}(0,I_{d_y})$ are independent standard normal random variables,$\beta_{t_n} \in \mb{R}^{d_{y}}$ is an observation-channel parameter (likelihood control) at time $t_n$ that  determines the measurement variance and may itself be chosen as a control variable with an associated cost, $A \in \mb{R}^{d_x\times d_x}, B \in  \mb{R}^{d_x\times d_{\alpha}},  \alpha_t \in \mathbb{R}^{d_\alpha},  \sigma\in\mathbb{R}^{d_x\times d_w},W_t \in \mathbb{R}^{d_w}, C \in \mb{R}^{d_y\times d_x}, m_0 \in  \mb{R}^{d_x}, \Sigma_0 \in  \mb{R}^{d_x\times d_x} $.
	The goal is to find the optimal control policy $\alpha^*,\beta^{*}$ that minimizes the expected cost functional:
	\begin{equation}
		\label{eq:POSOC_LQG_cost_functional}
		\begin{aligned}
			J(\alpha,\beta) &= \mb{E} \left[ \int_0^T  \frac{1}{2} \left( X_t^{\top} \, Q \, X_t + \alpha_t^{\top} \,R \,\alpha_t \right) \dl t \, + \, \frac{1}{2}  X_T^{\top} \, Q_T \, X_T +  \sum_{n=1}^{N_o} c_n\left(\beta_{t_n}\right)  \right],\\
		\end{aligned}
	\end{equation}
	where $Q \in \mb{R}^{d_x \times d_x}$, $Q_T \in \mb{R}^{d_x \times d_x}$, and $R \in \mb{R}^{d_{\alpha} \times d_{\alpha}}$ are positive semi-definite matrices.
	The discrete cost function $c_n(\beta_{t_n})$ can be defined as $$  c_n\left(\beta_{t_n}\right) = \sum_{i = 1}^{d_y} \frac{\kappa_{n,i}}{\beta_{t_n,i}} =\trace(\diag(\kappa_n) \diag^{-1}(\beta_{t_n}))$$
	where $\kappa \in \mb{R}^{d_{y}}$. It is intuitive that for smaller values of $\beta_{t_n}$  we will have  better observation of the state process $X_t$ and thus a better control policy $\alpha^*$. However, this will also lead to a larger cost associated with the observation process.
	This problem  will be referred  to as  \textbf{Partially Observed  Stochastic Optimal Control (POSOC-LQG)} throughout the paper. 
\end{exam}
In the above observation model, the sequence $(\xi_n)_n$ with
$\xi_n \sim \mcal N(0,I_{d_y})$ is the \emph{exogenous} measurement noise.
Its law is fixed and does not depend on the controller.
The process $\beta_{t_n} \in \R^{d_y}$ does not alter the distribution of
$\xi_n$ itself; instead, it determines how this noise is scaled into the
observation and thus parametrizes the conditional law (likelihood) of
$Y_{t_n}$ given $X_{t_n}$. Smaller values of $\beta_{t_n}$ correspond to
more informative measurements (lower observation variance), but are
penalized through the term $c_n(\beta_{t_n})$ in
\eqref{eq:POSOC_LQG_cost_functional}.
For this reason, it is more precise to view $\beta$ as a
\emph{likelihood control} (or observation-channel control) rather than as a
control of the driving noise: we do not control $(\xi_n)_n$, only the way
in which it enters the measurement process via the likelihood.

\subsection{Belief-state reformulation (fully observed on measures)}

For $t\in(t_n,t_{n+1})$, define the (pathwise) filtering distribution, also referred to as the 
\emph{belief} (or \emph{information state} in some control-theoretic literature) 
$$
\mu_t(\cdot\mid \mf Y^{[n]}) := \mcal L\!\left(X_t\mid \mf Y^{[n]}\right)\in\mcal P(\mb R^{d_x}).
$$
Here $\mathcal{L}(\cdot)$ denotes the  law of a random variable, and $\mathcal{P}(\mathbb{R}^{d_x})$ is the set of Borel probability measures on $\mathbb{R}^{d_x}$.

\begin{rem}[Filtrations vs.\ observation vectors]
	The observation filtration \eqref{eq:filtration-Y} satisfies $\mb F^Y_t=\mb F^Y_{t_n}=\sigma(\mf Y^{[n]}) $ for $t\in[t_n,t_{n+1})$, hence 
	$$
	\mu_t(\cdot\mid \mf Y^{[n]})
	=\mcal L\!\left(X_t\mid \mb F^Y_t\right).
	$$
	All developments below can therefore be phrased equivalently in filtration notation (with controls $\mb F^Y$–adapted). We retain the concrete path notation $\mf Y^{[n]}$ to make the jump updates at $\{t_n\}$ explicit and to highlight regression-style computations at observation times.
\end{rem}

We let $\mcal G_\alpha$ denote the generator of \eqref{eq:POSOC_dynamics} acting on test functions 
and $\mcal G_\alpha^*$ denote its adjoint acting on measures.  
We also use the pairing $\lrang{\varphi}{\nu}:=\int \varphi(x)\,\nu(\dl x)$.
Between observation times, the belief evolves deterministically:
\[
\dot{\mu}_t(\cdot\mid\mf Y^{[n]}) \;=\; \mcal G_{\alpha_t}^*\,\mu_t(\cdot\mid\mf Y^{[n]}),
\qquad t\in(t_n,t_{n+1}).
\]
At $t_{n}$, a Bayesian update maps the prior $\mu_{t_n^-}$ to the posterior $\mu_{t_n}$ via the likelihood $\pi_n$:
\[
\mu_{t_n}(\dl x\mid \mf Y^{[n-1]},y_n)
= \mcal K_{\beta_{t_n},y_n}\!\big(\cdot;\mu_{t_n^-}\big)
:= \frac{\pi_n\!\left(y_n\mid x,\mf Y^{[n-1]},\beta_{t_n}\right)}{L_n\!\left(y_n;\mf Y^{[n-1]},\beta_{t_n}\right)}\,\mu_{t_n^-}(\dl x),
\]
where the predictive normalizer is
\[
L_n\!\left(y;\mf Y^{[n-1]},\beta\right)
:= \int_{\mb R^{d_x}}\pi_n\!\left(y\mid x,\mf Y^{[n-1]},\beta\right)\,\mu_{t_n^-}(\dl x).
\]
Define the \emph{averaged} costs for any $\mu\in\mcal P(\mb R^{d_x})$:
\[
\tilde f_t(\mu,\alpha):=\lrang{f_t(\cdot,\alpha)}{\mu},\qquad
\tilde c_{t_n}(\beta,\mu):=\lrang{c_{t_n}(\cdot,\beta)}{\mu},\qquad
\tilde g(\mu):=\lrang{g(\cdot)}{\mu}.
\]
Taking expectation over observation paths, \eqref{eq:POSOC_cost_functional-X} becomes
\begin{align}
	\label{eq:POSOC_cost_functional_mu}
	J(u)
	= \mb E\!\left[
	\sum_{n=0}^N \int_{t_n}^{t_{n+1}}\!\tilde f_t\!\big(\mu_t(\cdot\mid \mf Y^{[n]}),\alpha_t\big)\,\dl t
	+ \sum_{n=1}^{N_o}\tilde c_{t_n}\!\big(\beta_{t_n},\mu_{t_n^-}(\cdot\mid \mf Y^{[n-1]})\big)
	+ \tilde g\big(\mu_T(\cdot\mid \mf Y^{[N]})\big)
	\right].
\end{align}
The corresponding HJB is a functional PDE on $\mcal P(\mb R^{d_x})$.  
In Gaussian settings, it reduces to a finite–dimensional HJB in the mean–covariance state.
The following result is adapted from \cite[Theorem~3.9]{bayer2024continuoustimestochasticoptimal},
with notation adjusted to our setting.

\noindent\textbf{Belief-space value functional.}
Fix $t\in(t_n,t_{n+1})$ and a belief $\mu\in\mathcal P(\mathbb R^{d_x})$.
Let $\mu_s^{u}$ denote the belief flow induced by a partially observed control
$u\in\mathcal U^{\mathrm{PO}}[t,T]$, with
$\mu_t^{u}=\mu$ and $k=\lfloor s\rfloor$ the unique index such that
$s\in[t_k,t_{k+1})$.

The policy-dependent belief cost function conditioned on $\mu$ is :
\begin{align}\label{eq:Vu-def-abstract}
	V_{t}\left(\mu; u\right)
	&:= \int_t^{t_{n+1}} \tilde f_s\big(\mu_s^{u},\alpha_s\big)\,\mathrm{d}s
	+ \tilde c_{t_{n+1}}\big(\beta_{t_{n+1}},\mu_{t_{n+1}^-}^{u}\big)\\
	&\quad+\mathbb E\Bigg[
	\int_{t_{n+1}}^T \tilde f_s\big(\mu_s^{u},\alpha_s\big)\,\mathrm{d}s
	+ \sum_{i = n+2}^{N_o}\tilde c_{t_i}\big(\beta_{t_i},\mu_{t_i^-}^{u}\big)+
	\tilde g\big(\mu_T^{u}\big)
	\;\Big|\; \mu_t^{u}=\mu
	\Bigg].
\end{align}

The optimal belief cost function conditioned on $\mu$ (value function\footnote{In the dynamic programming literature, the minimal expected cost as a function
	of the information state (here, the belief $\mu\in\mathcal P(\mathbb R^{d_x})$)
	is traditionally called the \emph{value function}. Strictly speaking,
	$V_t(\cdot)$ is a functional of the probability measure $\mu$, but we keep the
	standard terminology. }) is
\begin{equation}\label{eq:Vstar-def-abstract}
	V_t\left(\mu\right)
	:= \min_{u\in\mathcal U^{\mathrm{PO}}[t,T]} V_t(\mu ; u),
	\qquad t\in(t_n,t_{n+1}),\ \mu\in\mathcal P(\mathbb R^{d_x}).
\end{equation}

\begin{rem}[Belief and conditioning]
	There are two related viewpoints of the belief:
	
	(i) a generic belief $\mu\in\mathcal P(\mathbb R^{d_x})$; and
	
	(ii) a realized conditional law
	$\mu^{\mathbf y^{[n]}}:=\mathcal L(X_t\mid \mathbf Y^{[n]}=\mathbf y^{[n]})$
	for a fixed data realization $\mathbf y^{[n]}$.
	
The map $\Phi:\mathbf y^{[n]}\mapsto \mu^{\mathbf y^{[n]}}$ has image
$\Phi(\mathcal Y^{[n]})\subset \mathcal P(\mathbb R^{d_x})$, which is typically
a  submanifold of the full belief space.
\end{rem}
\begin{figure}[H]
	\centering
\begin{tikzpicture}[scale=2.5, >=stealth]
	
	\shade[ball color=gray!10, opacity=0.5] (0,0) ellipse (3 and 1.2);
	\draw (0,0) ellipse (3 and 1.2);
	\node at (2.0,1.3) {$\mathcal P(\mathbb R^{d_x})$};
	
	\coordinate (mu0)      at (-2,-0.3);
	\coordinate (mu1minus) at (-1.2,-0.5);
	\coordinate (mu1y1)    at (-0.15,0.4);
	\coordinate (mu2minus) at (0.8,0.3);
	\coordinate (mu2y1)    at (1.4,0.3);
	
	\draw[thick,->] (mu0) -- (mu1minus);
	
	\fill (mu0) circle (0.035);
	\node[below left=1pt] at (mu0) {$\mu_{t_0}^{\emptyset}$};
	
	\fill (mu1minus) circle (0.035);
	\node[below right=1pt] at (mu1minus) {$\mu_{t_1^-}$};
	
	\draw[thick,->,dashed] (mu1minus) -- (mu1y1);
	\fill (mu1y1) circle (0.025);
	\node[above right=1pt] at (mu1y1) {$\mu_{t_1}^{\mf y^{[1]}}$};
	
	\node[above=1pt] at (-0.8,-0.1)
	{$\mathcal K_{\beta,y_1}$};
	
	\coordinate (mu1y2) at (-0.3,-0.1);
	\coordinate (mu1y3) at (-0.4,-0.6);
	\coordinate (mu1y4) at (-0.2,-1.02);
	
	\coordinate (mu2minusy2) at (0.7,-0.1);
	\coordinate (mu2minusy3) at (0.7,-0.5);
	\coordinate (mu2minusy4) at (0.9,-1.02);
	
	\draw[thin,->,dotted] (mu1minus) -- (mu1y2);
	\draw[thin,->,dotted] (mu1minus) -- (mu1y3);
	\draw[thin,->,dotted] (mu1minus) -- (mu1y4);

	\draw[blue,thin,->,] (mu1y2)
	.. controls (0.0,-0.3) and (0.2,0.2) ..
	(mu2minusy2);
		\draw[blue,thin,->,] (mu1y3)
	.. controls (0.35,-0.65) and (0.2,-0.75) ..
	(mu2minusy3);
	\draw[blue,thin,->,] (mu1y4)
.. controls (0.35,-0.9) and (0.2,-0.75) ..
(mu2minusy4);
	
	\fill[gray] (mu2minusy2) circle (0.02);
	\fill[gray] (mu2minusy3) circle (0.02);
	\fill[gray] (mu2minusy4) circle (0.02);

		\draw[blue, thick]
	plot [smooth] coordinates {(mu2minusy4) (mu2minusy3) (mu2minusy2) (mu2minus)};

	\fill[gray] (mu1y2) circle (0.02);
	\fill[gray] (mu1y3) circle (0.02);
	\fill[gray] (mu1y4) circle (0.02);

	\node[right, gray] at (mu1y2) {$\mu_{t_1}^{y}$};
	\node[right, gray] at (mu1y3) {$\mu_{t_1}^{y'}$};
	\node[right, gray] at (mu1y4) {$\mu_{t_1}^{y''}$};
	
	\draw[blue, thick]
	plot [smooth] coordinates {(mu1y4) (mu1y3) (mu1y2) (mu1y1)};
	
	\node[blue, below right=1pt] at (0,.1)
	{$y \mapsto \mu_{t_1}^y$};
	
	\draw[thick,->] (mu1y1)
	.. controls (0.35,0.05) and (0.6,-0.05) ..
	(mu2minus);
	
	\fill (mu2minus) circle (0.035);
	\node[above left=1pt] at (mu2minus) {$\mu_{t_2^-}^{\mf y^{[1]}}$};
	
	\draw[thick,->,dashed] (mu2minus) -- (mu2y1);
	\fill (mu2y1) circle (0.025);
	\node[above right=1pt] at (mu2y1) {$\mu_{t_2}^{\mf y^{[2]}}$};
	
	\node[above=1pt] at (1.1,0.3)
	{$\mathcal K_{\beta,y_2}$};
	
	\coordinate (mu2y2) at (1.3,-0.05);
	\coordinate (mu2y3) at (1.35,-0.45);

	\coordinate (mu2y2y21) at (0.9,-0);
	\coordinate (mu2y2y22) at (0.87,-0.1);
	\coordinate (mu2y2y23) at (0.9,-0.2);

	\draw[thin,->,dotted] (mu2minusy2) -- (mu2y2y21);
	\draw[thin,->,dotted] (mu2minusy2) -- (mu2y2y22);
	\draw[thin,->,dotted] (mu2minusy2) -- (mu2y2y23);
	
	\node[left, gray] at (mu2minusy2) {$\mu_{t_2^-}^{y}$};
	\node[left, gray] at (mu2minusy3) {$\mu_{t_2^-}^{y'}$};
	\node[left, gray] at (mu2minusy4) {$\mu_{t_2^-}^{y''}$};
	
	\draw[ForestGreen, thick]
	plot [smooth] coordinates {(mu2y2y21) (mu2y2y22) (mu2y2y23)};
	
	\coordinate (mu2y2y31) at (0.85,-0.35);
	\coordinate (mu2y2y32) at (0.85,-0.45);
	\coordinate (mu2y2y33) at (0.9,-0.55);
	
	\draw[thin,->,dotted] (mu2minusy3) -- (mu2y2y31);
	\draw[thin,->,dotted] (mu2minusy3) -- (mu2y2y32);
	\draw[thin,->,dotted] (mu2minusy3) -- (mu2y2y33);
		\draw[ForestGreen, thick]
	plot [smooth] coordinates {(mu2y2y31) (mu2y2y32) (mu2y2y33)};

	\coordinate (mu2y2y41) at (1.05,-0.85);
	\coordinate (mu2y2y42) at (1.0,-0.95);
	\coordinate (mu2y2y43) at (1.0,-1.02);
	
	\draw[thin,->,dotted] (mu2minusy4) -- (mu2y2y41);
	\draw[thin,->,dotted] (mu2minusy4) -- (mu2y2y42);
	\draw[thin,->,dotted] (mu2minusy4) -- (mu2y2y43);
	
		\draw[ForestGreen, thick]
	plot [smooth] coordinates {(mu2y2y41) (mu2y2y42) (mu2y2y43)};
	
	\draw[thin,->,dotted] (mu2minus) -- (mu2y2);
	\draw[thin,->,dotted] (mu2minus) -- (mu2y3);
	
	\fill[gray] (mu2y2) circle (0.02);
	\fill[gray] (mu2y3) circle (0.02);
	
	\node[right, gray] at (mu2y2) {$\mu_{t_2}^{[y_1,y_2']}$};
	\node[right, gray] at (mu2y3) {$\mu_{t_2}^{[y_1,y_2'']}$};
	
	\draw[red, thick]
	plot [smooth] coordinates {(mu2y3) (mu2y2) (mu2y1)};
	
	\node[red, above left=1pt] at (2.55,-0.25)
	{$y \mapsto \mu_{t_2}^{[y_1,y]}$};

	\node at (1.55,-0.85) {\scriptsize $\mf y^{[1]} = (y_1)$};
	\node at (1.55,-1.05) {\scriptsize $\mf y^{[2]} = (y_1,y_2)$};
	
\begin{scope}
	\fill[ForestGreen!40, opacity=0.25]
	plot [smooth cycle] coordinates {
		(mu2minusy4)  
		(mu2minusy3)
		(mu2minusy2)
		(mu2minus)    
		(mu2y1)       
		(mu2y2)
		(mu2y3)       
	};
\end{scope}
\node[ForestGreen!60!black] at (1.5,-0.65)
{\scriptsize $(y_1,y_2)\mapsto \mu_{t_2}^{[y_1,y_2]}$};
	
\end{tikzpicture}

\caption{Belief dynamics on a data–indexed submanifold of $\mathcal P(\mathbb R^{d_x})$.
	Starting from an initial belief $\mu_{t_0}^{\emptyset}$, the black curve shows the
	prediction flow between observation times, yielding $\mu_{t_1^-}$ and
	$\mu_{t_2^-}^{\mathbf y^{[1]}}$.
	At each observation time $t_i$, the dashed arrow labeled $\mathcal K_{\beta,y_i}$
	represents the Bayesian update for the realized data $\mathbf y^{[i]}$,
	producing the posteriors $\mu_{t_1}^{\mathbf y^{[1]}}$ and
	$\mu_{t_2}^{\mathbf y^{[2]}}$.
	The blue and red curves represent, respectively, the one–dimensional families
	$y_1 \mapsto \mu_{t_1}^{(y_1)}$ and $y_2 \mapsto \mu_{t_2}^{[y_1,y_2]}$,
	while the green shaded surface shows the two–dimensional data–indexed belief
	manifold $(y_1,y_2)\mapsto \mu_{t_2}^{[y_1,y_2]}$ embedded in
	$\mathcal P(\mathbb R^{d_x})$.}

	\label{fig:belief-manifold-trajectory}
\end{figure}
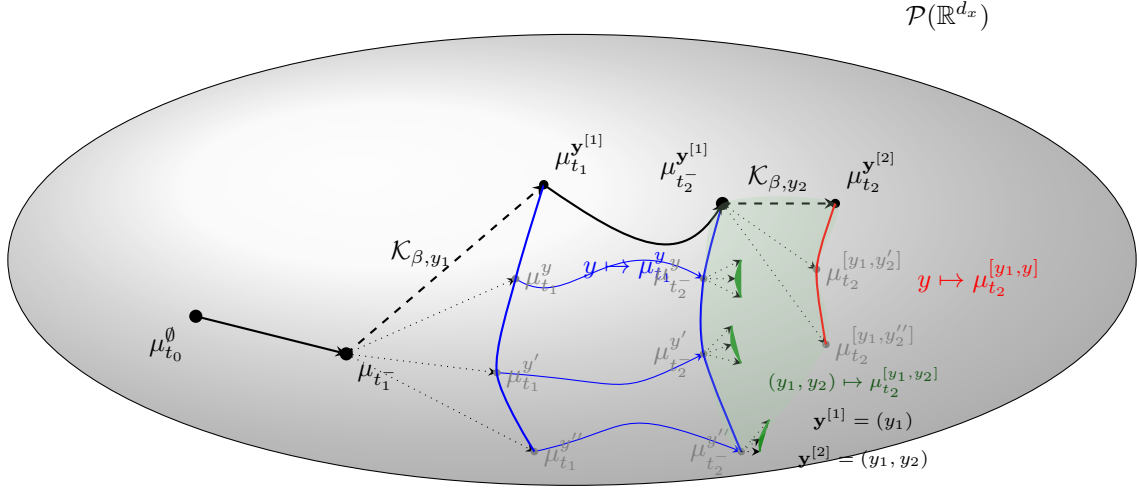

We recall the following definition and theorem from \cite{bayer2024continuoustimestochasticoptimal}, adapted to our notation.
\begin{defi}[Class $\mathcal{S}^{1,1}\left(\mathcal{P}\left(\mathbb{R}^d\right)\right)$ ]
 We say that a function $\Phi \in \mathcal{S}^{1,1}\left(\mathcal{P}\left(\mathbb{R}^d\right)\right)$ if there is a continuous version of the flat derivative $\frac{\delta \Phi}{\delta \mu}(\mu, x)$ such that
 \begin{itemize}
	\item  the mapping $(\mu, x) \mapsto \frac{\delta \Phi}{\delta \mu}(\mu, x)$ is jointly continuous w.r.t. $(\mu, x)$,
	\item  the mapping $x \mapsto \frac{\delta \Phi}{\delta \mu}(\mu, x)$ is twice continuously differentiable with bounded first and second order derivatives.
	 \end{itemize}
\end{defi}

\begin{rem}[Attainment of minima]
	Throughout we assume that all optimization problems under consideration admit minimizers.
	In particular, the infima in the next definitions are attained, so that
	we may write $\min$ and $\argmin$ \eqref{eq:optimal-belief-feedback}.
\end{rem}
\begin{theo}[Belief-space HJB]
	Let $V_t:\mathcal P(\mathbb R^{d_x})\to\mathbb R$ be value function defined in~\eqref{eq:Vstar-def-abstract}.
	Suppose that $V_t$ is differentiable w.r.t. the time variable $t$ and is in  $\mathcal{S}^{1,1}\left(\mathcal{P}\left(\mathbb{R}^d\right)\right)$ w.r.t. $\mu$. Then it satisfies the HJB equation, for all $n=0,\dots,N_o$ and all $t\in(t_n,t_{n+1})$, 
\begin{subequations}\label{eq:belief-HJB}
	\begin{align}
		\frac{\partial V_t}{\partial t}(\mu)
		&+\inf_{(\alpha,\beta)\in (\mcal A, \mcal B) }\Big\{
		\langle\mathcal G_\alpha \tfrac{\delta V_t}{\delta\mu}(\mu,\cdot),\,\mu\rangle
		+ \tilde f_t(\mu,\alpha)\Big\}
		= 0, \qquad t\in(t_n,t_{n+1}),\ n = 0,\dots,N, 
		\label{eq:belief-HJB-cont}\\[0.4em]
		V_{t_n^-}(\mu)
		&=\inf_{(\alpha,\beta)\in (\mcal A, \mcal B)}\left\{\tilde c_{t_n}(\beta,\mu)
		+ \mathbb E\!\left[V_{t_n}\!\big(\mathcal K_{\beta,Y_{t_n}}(\cdot;\mu)\big)\right]\right\},
		\qquad n = 1,\dots,N,
		\label{eq:belief-HJB-jump}\\[0.4em]
		V_T(\mu)
		&=\tilde g(\mu).
		\label{eq:belief-HJB-terminal}
	\end{align}
\end{subequations}

Fix, for each $t\in[0,T]$, a function 
$u^{*}_t : \mcal{P}\left(\mathbb{R}^d\right) \to (\mcal A, \mcal B)$
defined by
\begin{equation}\label{eq:optimal-belief-feedback}
	u^{*}_t(\mu) \in 
	\begin{cases}
		\displaystyle
		\argmin_{(\alpha,\beta)\in (\mcal A, \mcal B)}\Big\{
		\langle\mathcal G_\alpha \tfrac{\delta V_t}{\delta\mu}(\mu,\cdot),\,\mu\rangle
		+ \tilde f_t(\mu,\alpha)\Big\},
		& t\in(t_n,t_{n+1}),\ n = 0,\dots,N,\\[0.6em]
		\displaystyle
		\argmin_{(\alpha,\beta)\in (\mcal A, \mcal B)} \left\{\tilde c_{t_n}(\beta,\mu)
		+ \mathbb E\!\left[V_{t_n}\!\big(\mathcal K_{\beta,Y_{t_n}}(\cdot;\mu)\big)\right]\right\},
		& t=t_n,\ n = 1,\dots,N.
	\end{cases}
\end{equation}

provided  the set of minimizers in \eqref{eq:optimal-belief-feedback} is not empty.
\end{theo}

\subsection{Cost-to-go functions under different information structure}
For $t\in(t_n,t_{n+1})$ : 
\begin{itemize}
\item for a state $x\in\mb R^{d_x}$, and observation path $\mf y^{[n]}$, define for a fixed policy $u\in\mcal U^{\mathrm{R}}[t,T]$:
\begin{equation}\label{eq:UR-def}
	U_{t}^{\mathrm{R}}\left(x,\mf y^{[n]} ;\; u\right)
	:= \mb E\!\left[
	\int_t^T f_\tau(X_\tau,\alpha_\tau)\,\dl \tau
	+ \sum_{\{i:\,t_i\ge t\}} c_{t_i}\!\left(X_{t_i^-},\beta_{t_i}\right)
	+ g(X_T)\;\Big|\; (X_t,\mf Y^{[n]})=(x,\mf y^{[n]})
	\right].
\end{equation}

\item for a state $x\in\mb R^{d_x}$, define for a fixed policy $u\in\mcal U^{\mathrm{FO}}[t,T]$:
\begin{equation}\label{eq:UFO-def}
	U_t^{\mathrm{FO}}\left(x;\; u\right)
	:= \mb E\!\left[
	\int_t^T f_\tau(X_\tau,\alpha_\tau)\,\dl \tau
	+ \sum_{\{i:\,t_i\ge t\}} c_{t_i}\!\left(X_{t_i^-},\beta_{t_i}\right)
	+ g(X_T)\;\Big|\; X_t=x
	\right].
\end{equation}

\item For a state $x\in\mb R^{d_x}$, and observation path $\mf y^{[n]}$, define for a fixed policy $u\in\mcal U^{\mathrm{PO}}[t,T]$:
\begin{equation}\label{eq:UPO-def}
	U_{t}^{\mathrm{PO}}\left(\mf y^{[n]} ;\; u\right)
	:= \mb E\!\left[
	\int_t^T f_\tau(X_\tau,\alpha_\tau)\,\dl \tau
	+ \sum_{\{i:\,t_i\ge t\}} c_{t_i}\!\left(X_{t_i^-},\beta_{t_i}\right)
	+ g(X_T)\;\Big|\; \mf Y^{[n]}=\mf y^{[n]}
	\right].
\end{equation}
\end{itemize}

\begin{rem}[Interpreting  $\mathcal U^{\mathrm R}$ via randomized policies ]
	\label{rem:UR-seed}
	It is useful to view $\mathcal U^{\mathrm R}$ as 
	\emph{state–feedback with access to a random seed}.
	
	Under our assumptions, a controller in $\mathcal U^{\mathrm R}$ observes both the
	state history $X_{[0,t]}$ and the past measurements $Y_{t_1},\dots,Y_{t_{\nu(t)}}$.
	In the LQG example,
	\[
	Y_{t_n} \;=\; C X_{t_n} + \diag(\beta_{t_n})\,\xi_n,
	\]
	with i.i.d.\ Gaussian noise $(\xi_n)_n$ independent of $(X_0,W)$. For fixed
	$(X_{t_n},\beta_{t_n})$, each new measurement $Y_{t_n}$ carries new, independent random input 
	coming from $\xi_n$, and from this one can (via measurable transforms) generate
	uniform random variables, bits, etc.\ and hence implement any randomized decision
	rule.\footnote{In the
		linear-Gaussian case with $\diag(\beta_{t_n})$ invertible we can write
		$\xi_n = \diag(\beta_{t_n})^{-1}\big(Y_{t_n} - C X_{t_n}\big)$, so the extra
		information in $\mathcal F^{X,Y}_t$ beyond $\mathcal F^X_t$ is exactly the
		sequence $(\xi_1,\dots,\xi_{\nu(t)})$, and
		one can rewrite the joint filtration as
		$\mathcal F^{X,Y}_t = \mathcal F^X_t \vee \sigma(\xi_1,\dots,\xi_{\nu(t)})$. In more general observation models one
		usually assumes that the conditional law of $Y_{t_n}$ given $(X_{t_n},\beta_{t_n})$
		is non-degenerate; this is enough to treat $Y_{t_n}$ as a source of independent
		randomness for randomization.}
	
	A policy $u\in\mathcal U^{\mathrm R}$ can therefore:
	\begin{itemize}
		\item ignore this extra randomness and choose a \textbf{deterministic}
		state–feedback action based only on $X_{[0,t]}$ (recovering
		$\mathcal U^{\mathrm{FO}}$), or
		\item use the observation exogenous driving noise as a \textbf{random seed} to sample an
		action according to a state-dependent distribution.
	\end{itemize}
	In this sense, $\mathcal U^{\mathrm R}$ contains pure strategies as a special case and also allows mixed strategies
	through the observation exogenous driving  noise.
\end{rem}

\begin{enumerate}[(i)]
	\item \emph{Full-information benchmark value function} (state–feedback):
	\begin{equation}\label{eq:control-FO-def}
		\underline u_t(x)
		\in \argmin_{u\in\mcal U^{\mathrm{FO}}[t,T]} 	U_t^{\mathrm{FO}}(x;\; u),
	\end{equation}
	\begin{equation}\label{eq:value-FO-def}
		\underline U_t(x) := U_t^{\mathrm{FO}}(x;\; \underline u(x)).
	\end{equation}
	\item \emph{Partial-information value function.} 
	\begin{equation}\label{eq:control-PO-def}
		\tilde u_t \left( \mf y^{[n]}\right)\in\argmin_{u\in\mcal U^{\mathrm{PO}}[t,T]}
		U_{t}^{\mathrm{PO}}\left(\mf y^{[n]} ;\;  u \right),
	\end{equation}
	\begin{equation}\label{eq:value-PO-def}
		\tilde U _t \left( \mf y^{[n]}\right):= 
		U_{t}^{\mathrm{PO}}\left(\mf y^{[n]} ;\; \tilde u_t \left( \mf y^{[n]}\right) \right),
	\end{equation}
	and since $\mcal{U}^{\mathrm {PO}} \subseteq \mcal{U}^{ \mathrm R}$ we can define the auxiliary value induced by $\tilde u$ over $U_{t}^{[n],\mathrm{R}}$ : 
	\begin{equation}\label{eq:Urealized-def}
		\widehat U _t\left(x, \mf y^{[n]}\right)
		= U_t^{\mathrm {R}}\left(x,\mf y^{[n]} ;\; \tilde u_t \left(\mf y^{[n]}\right)   \right).
	\end{equation}
\end{enumerate}
The minimization is over \emph{policies} 
$u=(\alpha_s,\beta_{s})_{s\in[t,T]}$
(adapted on $[t,T]$ to the observation filtration), i.e., the entire control path from $t$ to $T$, 
not only the instantaneous action at time $t$.

\begin{propo}[Envelope inequality under randomized enlargement]
	\label{prop:envelope-ineq}
	Fix $t\in(t_n,t_{n+1})$, $x\in\R^{d_x}$ and an observation history $\mf y^{[n]}$. Then
	\begin{equation}\label{eq:U-ineq-envelope-claim}
		\underline U_t(x)\;\le\; \widehat U_t(x,\mf y^{[n]}),
	\end{equation}
	where $\underline U_t$ is the fully observed (state–feedback) value and $\widehat U_t^{[n]}$ is the realized pathwise value obtained by applying the partially observed optimizer $\tilde u(\mf y^{[n]})\in\mcal U^{\mathrm{PO}}[t,T]$ and evaluating it in the larger class $\mcal U^{\mathrm{R}}[t,T]$.
\end{propo}

\begin{proof}
	Let
	$$
	\check{U}_t(x,\mf y^{[n]}):=\min_{u\in\mcal U^{\mathrm R}[t,T]} U_t^{\mathrm R}(x,\mf y^{[n]};u).
	$$

	By the randomization approach to stochastic control, the value of the randomized problem equals that of the original fully observed problem :
	\begin{equation}\label{eq:FOeqR}
		\check{U}_t(x,\mf y^{[n]}) \;=\; 	\underline {U}_t(x)
	\end{equation}
	See, e.g., Fuhrman~\cite[Sec.~2.3 \& Sec.~5]{FuhrmanRandomization2025}, where the auxiliary (randomized) problem is constructed and the equality  is proved.
	Let $\tilde u(\mf y^{[n]})\in\mcal U^{\mathrm{PO}}[t,T]$ be a minimizer for the partially observed problem at history $\mf y^{[n]}$. Viewing $\tilde u(\mf y^{[n]})$ as an element of $\mcal U^{\mathrm R}$ that ignores $X$
	and because $\check{U}_t(x,\mf y^{[n]})$ is a minimum 
	$$
	\check{U}_t(x,\mf y^{[n]})
	\;\le\;
	U_t^{\mathrm R}\big(x,\mf y^{[n]};\,\tilde u_t(\mf y^{[n]})\big)
	=\widehat U_t(x,\mf y^{[n]}).
	$$
	Combining with \eqref{eq:FOeqR} yields \eqref{eq:U-ineq-envelope-claim}.
\end{proof}
For  realized data  $\mathbf y^{[n]}$ and  any fixed $u\in\mcal U^{\mathrm{PO}}[t,T]$ ,
\begin{equation}
	\label{eq:fixed-u-consistency}
	\begin{aligned}
		V_t(\mu^{\mathbf y^{[n]}};u) &= U_t^{ \mathrm{PO} }\left(\mf y^{[n]};u\right)=\left\langle U_t^{ \mathrm{R} }\left(\cdot,\mf y^{[n]};u\right),\,\mu^{\mf y^{[n]}}\right\rangle,
	\end{aligned}
\end{equation}
and, in particular, $V_t\left(\mu^{\mf y^{[n]}}\right)= \tilde U_t\left(\mf y^{[n]}\right),=\left\langle \widehat U_t\left(\cdot,\mf y^{[n]}\right),\mu^{\mf y^{[n]}}\right\rangle$.

Applying a needle variation to the auxiliary cost-to-go \eqref{eq:UR-def} yields a belief–weighted stationarity condition. 
	\begin{equation}\label{eq:ubar-mu}
		\overline u\left(\mf{y}^{[n]}\right) \in \begin{cases}
			\displaystyle
			\argmin_{\alpha\in\mcal A^{\mathrm{PO}}}
			\int_{\mb R^{d_x}}\!\Big(f_t(x,\alpha)+\mcal G_\alpha\, \overline U_t(x,\mf y^{[n]})\Big)\,
			\overline \mu_t(\dl x\mid \mf y^{[n]}), \qquad \forall t \in (t_n,t_{n+1})\\
			\displaystyle
			\argmin_{\beta\in\mcal B^{\mathrm{PO}}}\!
			\int_{\mb R^{d_x}}\!\Bigg(c_{t}(x,\beta)
			+ \int_{\mb R^{d_y}} \overline U_t\big(x,(\mf y^{[n]},y)\big)\,
			\pi_{n+1}(y\mid x,\mf y^{[n]},\beta)\,\dl y\Bigg)\,
			\overline \mu_{t^-}(\dl x\mid \mf y^{[n]}),\quad t=t_{n+1}
		\end{cases},
	\end{equation}
	Then, regardless of the global optimality, the $\mu$–average is an \emph{upper bound} on the optimal belief value:
	\begin{equation}\label{eq:V-ineq-envelope}
		V_t(\mu^{\mf y^{[n]}})
		=\left\langle \widehat U_t\left(\cdot,\mf y^{[n]}\right),\mu^{\mf y^{[n]}}\right\rangle
		\ \le\ 
		\left\langle \overline U_t(\cdot;\mu^{\mf y^{[n]}}),\,\mu^{\mf y^{[n]}}\right\rangle,
	\end{equation}
	where under the necesssary convexity  (the control set is convex and the running and terminal costs
	are convex in the state/control so that the optimal control is unique), these two quantities coincide.

	\subsection{Auxiliary optimality representation}
	\label{subsec:aux-optimality}
	
	We now derive first–order (needle–variation) necessary conditions for a locally optimal policy in the partially observed class. Fix an observation history $\mf y^{[n]}=(y_{t_1},\dots,y_{t_n})$ and consider the auxiliary cost-to-go :
	\[
	U_{t}^{\mathrm{R}}\left(x,\mf y^{[n]} ;\; u\right)
	:= \mb E\!\left[
	\int_t^T f_\tau(X_\tau,\alpha_\tau)\,\dl \tau
	+ \sum_{\{i:\,t_i\ge t\}} c_{t_i}\!\left(X_{t_i^-},\beta_{t_i}\right)
	+ g(X_T)\;\Big|\; (X_t,\mf Y^{[n]})=(x,\mf y^{[n]})
	\right],
	\]
	as in \eqref{eq:UR-def}, where $u=(\alpha,\beta)\in\mcal U^{\mathrm{PO}}[t,T]$.
	The following statement characterizes the stationarity of a locally optimal triple.
	
	\begin{propo}[Auxiliary optimality system representation]
		\label{prop:aux-opt}
		Let $\bigl(\bar u,\overline \mu,\overline U\bigr)$ be a policy–state–cost triple in the partially observed class, with $\bar u=(\bar\alpha,\bar\beta)$, such that $\bar u$ is locally optimal in the sense of needle variations around $(t,\mf y^{[n]})$,
		i.e., for every admissible policy $u\in\mcal U^{\mathrm{PO}}[t_n,T]$, every $\ve>0$ small, and every measurable
		neighbourhood $B_\ve^{(n)}\subset(\R^{d_y})^n$  and
		$B_\ve^{(n)}\downarrow\{\mf y^{[n]}\}$ as $\ve\downarrow0$, the spiked policy $u^\ve$ defined by
		\[
		u_s^\ve :=
		\begin{cases}
			u_s, & s\in[t,t+\ve)\ \text{and}\ \mf Y^{[n]}\in B_\ve^{(n)},\\
			\bar u_s, & \text{otherwise},
		\end{cases}
		\]
		(and analogously when the spike is applied only at a single observation time $t_k$)
		satisfies $J(u^\ve)\ge J(\bar u)$.
		
		Then, for each $n=0,\dots,N_o$ and all $t\in(t_n,t_{n+1})$, the pair
		$\bigl(\overline U_t(\cdot,\mf y^{[n]}),\,\overline \mu_t(\cdot\mid \mf y^{[n]})\bigr)$ solves 
		\begin{equation}
			\label{eq:aux_system}
			\begin{cases}
				\overline{\mu}_0(x) = \mu_0(x),\\[1.1ex]
				\dot { \overline{\mu}}_t (\cdot\mid \mf y^{[n]})
				= \mcal G_{\bar\alpha_t}^{*}\,\overline{\mu}_t(\cdot\mid \mf y^{[n]}),& t\in(t_n,t_{n+1}),\\[0.9ex]
				\overline{\mu}_{t_{n+1}}(x\mid \mf y^{[n]},y)
				= \dfrac{\pi_{n+1}(y\mid x,\mf y^{[n]},\bar\beta_{n+1})}{L_{n+1}(y;\mf y^{[n]},\bar\beta_{n+1})}\,
				\overline{\mu}_{t_{n+1}^-}(x\mid \mf y^{[n]}),& n=0,\dots,N_o-1,\\[1.1ex]
				\overline U_T(x,\mf y^{[N]}) = g(x),\\[0.6ex]
				\dot {\overline U}_t(x,\mf y^{[n]})
				= -\Big(\mcal G_{\bar\alpha_t}\,\overline U_t(x,\mf y^{[n]}) + f_t(x,\bar\alpha_t)\Big),
				& t\in(t_n,t_{n+1}),\\[0.9ex]
				\displaystyle \overline U_{t_{n+1}^-}(x,\mf y^{[n]})
				= c_{t_{n+1}}(x,\bar\beta_{n+1})
				+ \int_{\mb R^{d_y}} \overline U_{t_{n+1}}\!\big(x,(\mf y^{[n]},y)\big)\,
				\pi_{n+1}(y\mid x,\mf y^{[n]},\bar\beta_{n+1})\,\dl y,& n=0,\dots,N_o-1,\\[1.1ex]
				\displaystyle
				\bar\alpha_t \in \argmin_{\alpha\in\mcal A}
				\int_{\mb R^{d_x}}\!\Big(f_t(x,\alpha)+\mcal G_\alpha\,\overline U_t(x,\mf y^{[n]})\Big)\,
				\overline \mu_t(\dl x\mid \mf y^{[n]}),& t\in(t_n,t_{n+1}),\\[1.1ex]
				\displaystyle
				\bar\beta_{n+1}\in\argmin_{\beta\in\mcal B}\!
				\int_{\mb R^{d_x}}\!\Bigg(c_{t_{n+1}}(x,\beta)
				+ \int_{\mb R^{d_y}} \overline U_{t_{n+1}}\!\big(x,(\mf y^{[n]},y)\big)\,
				\pi_{n+1}(y\mid x,\mf y^{[n]},\beta)\,\dl y\Bigg)\,
				\overline{\mu}_{t_{n+1}^-}(\dl x\mid \mf y^{[n]})& n=0,\dots,N_o-1.
			\end{cases}
		\end{equation}
		\noindent
		In words: between observation dates, $\overline U$ evolves backward by the Kolmogorov backward operator associated with the controlled diffusion, the belief $\overline\mu$ evolves forward by the adjoint (Fokker–Planck) operator, the continuous-time control $\bar\alpha_t$ minimizes the \emph{belief-averaged} Hamiltonian, and at each observation date $t_{n+1}$ the pair $(\overline U,\overline\mu)$ undergoes a Bayesian jump driven by the observation kernel, while $\bar\beta_{n+1}$ minimizes the belief-averaged one-step cost.
	\end{propo}
	
	\noindent\emph{Proof.} See Appendix~\ref{app:proof-aux-opt}.
	
	\subsection{Lagrangian optimality representation}
	\label{subsec:POSOC_KKT}
	We regard the belief–space problem
	\[
	\min_{u\in\mcal U^{\mathrm{PO}}} J(u)
	\]
	as a constrained optimisation over $(u,\mu)$, where $\mu$ is the conditional law of $X_t$ given the observation history. The constraints are:
	(i) between two observation dates the belief satisfies the forward equation driven by the generator with the chosen continuous control, and
	(ii) at observation dates the belief jumps according to Bayes’ rule with the discrete sensing/action variable.
	
	We relax these constraints by introducing an adjoint (costate) :
	\[
\lambda_t(x,\mf Y^{[n]}),\qquad \forall x\in\R^{d_x},
\quad\text{(measurable w.r.t.\ $\F_{t_n}^Y=\sigma(\mf Y^{[n]})$).}
 	\]
	and define the Lagrangian
	\begin{align*}
		\mcal L(u,\mu,\lambda)
		&= J(u)
		\;-\; \mb E\!\left[\sum_{n=1}^{N_o+1}\int_{t_{n-1}}^{t_n}
		\left\langle \lambda_t(\cdot,\mf Y^{[n]}),\,
		\dot\mu_t(\cdot\mid \mf Y^{[n]})-\mcal G_{\alpha_t}^*\,\mu_t(\cdot\mid \mf Y^{[n]})
		\right\rangle \dl t\right] \\
		&\quad -\; \mb E\!\left[\sum_{n=1}^{N_o}
		\left\langle \lambda_{t_n}(\cdot,\mf Y^{[n]}),\,
		\mu_{t_n}(\cdot\mid \mf Y^{[n]})-\mcal K_{\beta_{t_n},Y_{t_n}}\!\big(\cdot;\,\mu_{t_n^-}(\cdot\mid \mf Y^{[n-1]})\big)
		\right\rangle \right].
	\end{align*}
	
	We can express the resulting forward–backward structure in the following statement.
	
	\begin{propo}[Lagrangian optimality system representation]
		\label{prop:FB-belief-cases}
		Let $(\bar u,\bar\mu)$ be an optimal pair, with $\bar u=(\bar\alpha,\bar\beta)$, and let $\bar\mu$ be the corresponding conditional laws. Then there exists an adjoint family
		\[
		\lambda_t(\cdot,\mf y^{[n]}), \qquad t\in(t_n,t_{n+1}),\; n=0,\dots,N_o,
		\]
		such that, for any $\mf y^{[n]}$ with positive density, the following pathwise system holds:
		\begin{equation}
			\label{eq:FOC-belief}
			\begin{cases}
				\overline{\mu}_0(x) = \mu_0(x),\\[1.1ex]
				\dot{\bar\mu}_t(\cdot \mid \mf y^{[n]})
				= \mcal G_{\bar\alpha_t}^*\,\bar\mu_t(\cdot \mid \mf y^{[n]}),
				& t\in(t_n,t_{n+1}),\\[0.8ex]
				\displaystyle
				\bar\mu_{t_{n+1}}^{[n+1]}(x \mid \mf y^{[n]},y)
				= \dfrac{\pi_{n+1}(y \mid x,\mf y^{[n]},\bar\beta_{n+1})}
				{\displaystyle \int_{\R^{d_x}} \pi_{n+1}(y \mid x',\mf y^{[n]},\bar\beta_{n+1})
					\,\bar\mu_{t_{n+1}^-}(\dl x' \mid \mf y^{[n]})}
				\;\bar\mu_{t_{n+1}^-}(x \mid \mf y^{[n]}),
				& n=0,\dots,N_o-1,\\[1.4ex]
				\lambda_T(x,\mf y^{[N]}) = g(x),
				& \\[0.8ex]
				\dot\lambda_t(x,\mf y^{[n]})
				= -\big(\mcal G_{\bar\alpha_t}\lambda_t(x,\mf y^{[n]})
				+ f_t(x,\bar\alpha_t)\big),
				& t\in(t_n,t_{n+1}),\\[0.8ex]
				\displaystyle
				\lambda_{t_{n+1}^-}(x,\mf y^{[n]})
				= c_{t_{n+1}}\!\left(x,\bar\beta_{n+1}\right)
				+ \int_{\R^{d_y}}\!\lambda_{t_{n+1}}(x,(\mf y^{[n]},y))\,
				\pi_{n+1}(y \mid x,\mf y^{[n]},\bar\beta_{n+1})\,\dl y \\[1.0ex]
				\qquad\displaystyle
				- \int_{\R^{d_x}}\!\bar\mu_{t_{n+1}^-}(\dl x' \mid \mf y^{[n]})
				\int_{\R^{d_y}}\!\lambda_{t_{n+1}}(x',(\mf y^{[n]},y))\,
				\frac{\pi_{n+1}(y \mid x',\mf y^{[n]},\bar\beta_{n+1})\,
					\pi_{n+1}(y \mid x,\mf y^{[n]},\bar\beta_{n+1})}
				{\displaystyle \int_{\R^{d_x}}\pi_{n+1}(y \mid z,\mf y^{[n]},\bar\beta_{n+1})
					\,\bar\mu_{t_{n+1}^-}(\dl z \mid \mf y^{[n]})}\,\dl y,
				& n=0,\dots,N_o-1\\
				{\displaystyle \bar\alpha_t \in \argmin_{\alpha\in\mcal A}
					\int_{\mb R^{d_x}}\!\Big(f_t(x,\alpha)+\mcal G_\alpha\,\overline \lambda_t(x,\mf y^{[n]})\Big)\,
					\overline \mu_t(\dl x\mid \mf y^{[n]})},& t\in(t_n,t_{n+1}),\\[1.1ex]
				\displaystyle
				\bar\beta_{n+1}\in\argmin_{\beta\in\mcal B}\!
				\int_{\mb R^{d_x}}\!\Bigg(c_{t_{n+1}}(x,\beta)
				+ \int_{\mb R^{d_y}} \overline \lambda_{t_{n+1}}\!\big(x,(\mf y^{[n]},y)\big)\,
				\pi_{n+1}(y\mid x,\mf y^{[n]},\beta)\,\dl y\Bigg)\,
				\overline{\mu}_{t_{n+1}^-}(\dl x\mid \mf y^{[n]})& n=0,\dots,N_o-1.
			\end{cases}
		\end{equation}
		
	\end{propo}
	
	\noindent
	Variations in $\lambda$ recover the belief flow and Bayesian jump (first two lines); variations in $\mu$ give the backward adjoint equation and the adjoint jump (last three lines). The detailed derivation, including the localization in $(t,\mf y^{[n]})$ and the integration–by–parts step on each $(t_n,t_{n+1})$, is given in Appendix~\ref{app:proof-FB-belief}.

	Subsections~\ref{subsec:aux-optimality} and~\ref{subsec:POSOC_KKT} yield two
	optimality systems that look slightly different. In the \emph{auxiliary}
	system (Proposition~\ref{prop:aux-opt}), the backward variable is the
	conditional cost-to-go $U$ and the jump at observation times is exactly the
	Bayesian update of the value, averaged with the observation  likelihood. In the
	\emph{Lagrangian/KKT} system (Proposition~\ref{prop:FB-belief-cases}), the
	backward variable is the multiplier $\lambda$ attached to the belief
	dynamics, and at observation times an additional term appears in the jump
	equation. This term is not an artefact: it is the contribution of
	differentiating the \emph{normalising denominator} in Bayes’ rule, i.e.\ the
	fact that the posterior measure depends on the prior through a ratio.
	
	Despite this difference in the backward equations, the two systems induce the
	\emph{same} optimal controls and belief flow. In both formulations, the
	forward components $\bar\mu_t(\cdot\mid\mf y^{[n]})$ and the
	controls $(\bar\alpha_t,\bar\beta_{n+1})$ are exactly those of the original
	partially observed problem. The continuous-time control is characterised by
	\[
	\bar\alpha_t \in \argmin_{\alpha\in\mcal A}
	\int_{\R^{d_x}}\!\Big(f_t(x,\alpha)+\mcal G_\alpha\,\overline U_t(x,\mf y^{[n]})\Big)\,
	\overline \mu_t(\dl x\mid \mf y^{[n]})
	\]
	in the auxiliary formulation, and by
	\[
	\bar\alpha_t \in \argmin_{\alpha\in\mcal A}
	\int_{\R^{d_x}}\!\Big(f_t(x,\alpha)+\mcal G_\alpha\,\lambda_t(x,\mf y^{[n]})\Big)\,
	\overline \mu_t(\dl x\mid \mf y^{[n]})
	\]
	in the Lagrangian formulation. At each observation time $t_{n+1}$, the
	discrete control $\bar\beta_{n+1}$ similarly minimises the same
	belief–averaged one-step cost, with $\overline U_{t_{n+1}}$ replaced by
	$\lambda_{t_{n+1}}$.
	
	The extra term in the jump of $\lambda$ in
	\eqref{eq:FOC-belief} comes from the dependence of the normalising factor in
	Bayes’ rule on the prior belief $\bar\mu_{t_{n+1}^-}$. This correction
	modifies $\lambda$, but does not change the belief–averaged quantities that
	enter the minimisation problems in $\alpha$ and $\beta$. Once $\overline U$
	and $\lambda$ are linked through their respective backward equations, the
	functionals being minimised in both formulations coincide, and both systems
	therefore correspond to the same optimal controls $(\bar\alpha,\bar\beta)$
	and the same optimal belief trajectory $\bar\mu$.

\begin{rem}[Analogy with Zakai and Kushner-Stratonovich equations]
	Between observation times, the auxiliary system of Proposition~\ref{prop:aux-opt}
	propagates  $\bar\mu$ linearly, and at each observation time
	$t_{n+1}$, the belief is updated by weighting the prior $\bar\mu_{t_{n+1}^-}$ with
	the likelihood and then normalising. This \enquote{predict-then-normalise} structure
	is reminiscent of the continuous-time nonlinear filtering equations: the
	unnormalised Zakai equation, which evolves an unnormalised conditional density
	linearly, and the normalised Kushner-Stratonovich equation, which adds a
	nonlinear normalisation term; see, for example,
	\cite{BainCrisan2009}.
	By contrast, the Lagrangian/KKT system of
	Proposition~\ref{prop:FB-belief-cases} is written directly in terms of the
	normalised belief, and the jump of the adjoint carries an extra term coming from
	the Bayes normalising constant, in analogy with the nonlinear correction term
	in the Kushner-Stratonovich equation. We emphasise that the analogy is purely
	structural, since here observations are discrete in time: our updates are exact
	Bayes-rule jumps, not continuous stochastic correction terms. Thus, one should not literally interpret our adjoint system as a Zakai or Kushner--Stratonovich equation; it only mirrors the two-step \emph{propagate-and-update} pattern of nonlinear filtering.
	
\end{rem}


\subsection{Pontryagin formulation }
\label{subsec:Pontryagin}
For $U\in C_b^2(\mb R^{d_x})$, $\mu\in\mcal P(\mb R^{d_x})$
and $t\in(t_n,t_{n+1})$, define the continuous Hamiltonian
\begin{equation}\label{eq:Hc-def}
	\mcal H_t^{\mathrm c}(\mu,U)
	\;:=\; \inf_{\alpha\in\mcal A}\ \lrang{\,f_t(\cdot,\alpha)+\mcal G_\alpha U\,}{\,\mu\,}.
\end{equation}
At $t_{n+1}$, for $\mu^-:=\mu_{t_{n+1}^-}(\cdot\mid \mf y^{[n]})$ and any $\beta\in\mcal B$, we define the (pre-posterior) discrete Hamiltonian
\begin{equation}
	\label{eq:Hd-pre-def}
	\mcal H_{n+1}^{\mathrm d}(\mu^-,U^+)
	:= \inf_{\beta\in\mcal B}\Bigg(
	\lrang{c_{t_{n+1}}(\cdot,\beta)}{\mu^-}
	+\int_{\mb R^{d_y}} \!\lrang{\,U^+(\cdot,(\mf y^{[n]},y))\,}{\,\pi_{n+1}(y\mid \cdot,\mf y^{[n]},\beta)\,\mu^-\,}\,\dl y
	\Bigg).
\end{equation}

Here we denote the adjoint by $U$ to match the auxiliary pathwise system
\eqref{eq:aux_system} rewritten in Hamiltonian form:
\begin{equation}\label{eq:Ham-system-joint}
				\begin{cases}
	\overline{\mu}_0(x) = \mu_0(x),\\[1.1ex]
	\dot { \overline{\mu}}_t (x\mid \mf y^{[n]})
	=\dfrac{\delta \mcal H_t^{\mathrm c}}{\delta U}\left(\overline{\mu}_t\left(\cdot\mid\mf{y}^{[n]}\right),\overline{U}_t\left(\cdot,\mf{y}^{[n]}\right), x\right),& t\in(t_n,t_{n+1}),\\[0.9ex]
	\overline{\mu}_{t_{n+1}}(x\mid \mf y^{[n]},y)
	= \frac{\dfrac{\delta \mcal H_{n+1}^{\mathrm d}}{\delta U}\left( \displaystyle \overline{\mu}_{t_{n+1}^-}\left(\cdot\mid\mf{y}^{[n]}\right), \overline{U}_{t_{n+1}}\left(\cdot,\mf{y}^{[n]},y\right), x \right)}{\displaystyle \left(\int_{\mb{R}^{d_x}}  \dfrac{\delta \mcal H_{n+1}^{\mathrm d}}{\delta U}  \left(\overline{\mu}_{t_{n+1}^-}\left(\cdot\mid\mf{y}^{[n]}\right), \overline{U}_{t_{n+1}}\left(\cdot,\mf{y}^{[n]},y\right), x' \right) \dl x'   \right)} & n=0,\dots,N_o-1,\\[1.1ex]
	\overline U_T(x,\mf y^{[N]}) = g(x),\\[0.6ex]
	\dot {\overline U}_t(x,\mf y^{[n]})
	= -\dfrac{\delta \mcal H_t^{\mathrm c}}{\delta \mu}\left(\overline{\mu}_t\left(\cdot\mid\mf{y}^{[n]}\right),\overline{U}_t\left(\cdot,\mf{y}^{[n]}\right), x\right),
	& t\in(t_n,t_{n+1}),\\[0.9ex]
	\displaystyle \overline U_{t_{n+1}^-}(x,\mf y^{[n]})
	= \dfrac{\delta \mcal H_{n+1}^{\mathrm d}}{\delta \mu}\left( \displaystyle \overline{\mu}_{t_{n+1}^-}\left(\cdot\mid\mf{y}^{[n]}\right), \overline{U}_{t_{n+1}}\left(\cdot,\mf{y}^{[n]},y\right), x \right) & n=0,\dots,N_o-1,\\[1.1ex]
	\displaystyle
	\bar\alpha_t \in \argmin_{\alpha\in\mcal A}
	\int_{\mb R^{d_x}}\!\Big(f_t(x,\alpha)+\mcal G_\alpha\,\overline U_t(x,\mf y^{[n]})\Big)\,
	\overline \mu_t(\dl x\mid \mf y^{[n]}),& t\in(t_n,t_{n+1}),\\[1.1ex]
	\displaystyle
	\bar\beta_{n+1}\in\argmin_{\beta\in\mcal B}\!
	\int_{\mb R^{d_x}}\!\Bigg(c_{t_{n+1}}(x,\beta)
	+ \int_{\mb R^{d_y}} \overline U_{t_{n+1}}\!\big(x,(\mf y^{[n]},y)\big)\,
	\pi_{n+1}(y\mid x,\mf y^{[n]},\beta)\,\dl y\Bigg)\,
	\overline{\mu}_{t_{n+1}^-}(\dl x\mid \mf y^{[n]})& n=0,\dots,N_o-1.
\end{cases}
\end{equation}

This system is the natural Pontryagin forward-backward system on the
\emph{space of beliefs}. The state variable is the conditional law
$\overline\mu_t(\cdot\mid\mf y^{[n]})\in\mathcal P(\R^{d_x})$, and the
costate is the auxiliary cost-to-go $\overline U_t(\cdot,\mf y^{[n]})$.

For $t\in(t_n,t_{n+1})$, the continuous Hamiltonian
$\mathcal H_t^{\mathrm c}$ plays the usual role: its variational
derivatives with respect to the second and first argument generate,
respectively, the forward Fokker-Planck equation and the backward
Kolmogorov equation in the sense that, for any minimiser
$\bar\alpha_t\in\arg\min_{\alpha\in\mathcal A}
\langle f_t(\cdot,\alpha)+\mathcal G_\alpha U_t,\mu_t\rangle$,
$$
\dfrac{\delta \mcal H_t^{\mathrm c}}{\delta U}\left(\overline{\mu}_t\left(\cdot\mid\mf{y}^{[n]}\right),\overline{U}_t\left(\cdot,\mf{y}^{[n]}\right), x\right)
= \mathcal G_{\bar\alpha_t}^{*}\overline{\mu}_t\left(x\mid\mf{y}^{[n]}\right),
$$
$$
\dfrac{\delta \mcal H_t^{\mathrm c}}{\delta \mu}\left(\overline{\mu}_t\left(\cdot\mid\mf{y}^{[n]}\right),\overline{U}_t\left(\cdot,\mf{y}^{[n]}\right), x\right)
= \mcal G_{\bar\alpha_t}\,\overline U_t(x,\mf y^{[n]}) + f_t(x,\bar\alpha_t).
$$

At any observation time $t_{n+1}$, the discrete Hamiltonian
$\mathcal H_{n+1}^{\mathrm d}$ plays the exact analogue of this role:
its variational derivative with respect to $U$ gives the
\emph{unnormalised} posterior measure, and its derivative with respect
to $\mu$ gives the backward jump of the costate.

Then, at a minimiser $\bar\beta_{n+1}$ and for every $x\in\R^{d_x}$,
\[
\dfrac{\delta \mcal H_{n+1}^{\mathrm d}}{\delta U}\left( \displaystyle \overline{\mu}_{t_{n+1}^-}\left(\cdot\mid\mf{y}^{[n]}\right), \overline{U}_{t_{n+1}}\left(\cdot,\mf{y}^{[n]},y\right), x \right)= \pi_{n+1}(y \mid x,\mf y^{[n]},\bar\beta_{n+1})
\,\bar\mu_{t_{n+1}^-}(x \mid \mf y^{[n]}),
\]
so that the posterior belief is recovered by normalisation,
and
\[
\dfrac{\delta \mcal H_{n+1}^{\mathrm d}}{\delta \mu}\left( \displaystyle \overline{\mu}_{t_{n+1}^-}\left(\cdot\mid\mf{y}^{[n]}\right), \overline{U}_{t_{n+1}}\left(\cdot,\mf{y}^{[n]},y\right), x \right)
=c_{t_{n+1}}(x,\bar\beta_{n+1})
+ \int_{\mb R^{d_y}} \overline U_{t_{n+1}}\!\big(x,(\mf y^{[n]},y)\big)\,
\pi_{n+1}(y\mid x,\mf y^{[n]},\bar\beta_{n+1})\,\dl y,
\]
which is exactly the jump condition for $U$ in
\eqref{eq:Ham-system-joint}. In this sense, both the continuous-time
evolution and the discrete Bayesian updates are generated by the
\emph{same} Hamiltonian objects $(\mathcal H_t^{\mathrm c},
\mathcal H_{n+1}^{\mathrm d})$ on the belief space, just as in the
classical Pontryagin principle with a finite-dimensional state.

\medskip

\subsection{Relationship between the belief state value function and the adjoint variable}

Throughout this subsection we work under the convexity conditions discussed above, which ensure that the needle-optimal triple of Proposition~\ref{prop:aux-opt} is in fact globally optimal for the partially observed control problem. We denote this triple by $(\bar u,\bar\mu,\bar U)$, and we write
\[
\bar U_t(x,\mf y^{[n]})
\]
for the associated adjoint (cost-to-go) process given by Proposition~\ref{prop:aux-opt}. For consistency with the pathwise representation in~\eqref{eq:fixed-u-consistency}, we recall that, along the optimal triple, $\bar U_t(\cdot,\mf y^{[n]})$ coincides $\bar\mu_t^{\mf y^{[n]}}$-a.e.\ with the realized auxiliary value $U_t^{\mathrm R}(\cdot,\mf y^{[n]};\bar u)$.

The belief-space formulation admits two complementary viewpoints:

\begin{itemize}
	\item the \emph{dynamic programming} viewpoint, in which the value functional $V_t(\mu)$ solves the HJB equation~\eqref{eq:belief-HJB} on $\mcal P(\R^{d_x})$; and
	\item the \emph{Pontryagin} viewpoint, in which an optimal policy $\bar u$ is characterized by the auxiliary forward-backward system of Proposition~\ref{prop:aux-opt}, with adjoint $\bar U_t(x,\mf y^{[n]})$.
\end{itemize}

These two objects are linked in a way that is directly analogous to the finite-dimensional identity
$p_t = \nabla_x V(t,X_t^*)$, but \emph{only along the optimal belief path}.

\medskip

Fix $t\in(t_n,t_{n+1})$ and a realized observation history $\mf y^{[n]}$. For any fixed partially observed policy $u\in\mcal U^{\mathrm{PO}}[t,T]$, the consistency relation
\eqref{eq:fixed-u-consistency} reads
\[
V_t(\mu^{\mf y^{[n]}};u)
= U_t^{\mathrm{PO}}(\mf y^{[n]};u)
= \big\langle U_t^{\mathrm R}(\cdot,\mf y^{[n]};u),\,\mu^{\mf y^{[n]}}\big\rangle.
\]
In particular, for the optimal control $\bar u$ of Proposition~\ref{prop:aux-opt}, we denote by $\bar\mu^{\mf y^{[n]}}$ the associated optimal belief process and we have
\[
V_t\big(\bar\mu_t^{\mf y^{[n]}}\big)
= V_t\big(\bar\mu_t^{\mf y^{[n]}};\bar u\big)
= \big\langle U_t^{\mathrm R}(\cdot,\mf y^{[n]};\bar u),\,\bar\mu_t^{\mf y^{[n]}}\big\rangle.
\]
Since $\bar U_t(\cdot,\mf y^{[n]})$ and $U_t^{\mathrm R}(\cdot,\mf y^{[n]};\bar u)$ coincide $\bar\mu_t^{\mf y^{[n]}}$-a.e., this identity may also be written as
\begin{equation}\label{eq:V-envelope-again}
	V_t\big(\bar\mu_t^{\mf y^{[n]}}\big)
	= \big\langle \bar U_t(\cdot,\mf y^{[n]}),\,\bar\mu_t^{\mf y^{[n]}}\big\rangle,
	\qquad t\in(t_n,t_{n+1}).
\end{equation}

\begin{lem}[Envelope inequality at the optimal belief (minimization case)]
	\label{lem:envelope-ineq}
	Fix $(t,\mf y^{[n]})$ and consider the partially observed control problem started from $(t,\mf y^{[n]})$, with value functional
	\[
	V_t(\mu)
	:= \inf_{u\in\mcal U^{\mathrm{PO}}[t,T]} V_t(\mu;u),
	\qquad
	V_t(\mu^{\mf y^{[n]}};u)
	= \big\langle U_t^{\mathrm R}(\cdot,\mf y^{[n]};u),\,\mu^{\mf y^{[n]}}\big\rangle.
	\]
	Let $(\bar u,\bar\mu,\bar U)$ be the globally optimal triple of Proposition~\ref{prop:aux-opt}, and write $\bar\mu_t^{\mf y^{[n]}}$ for the corresponding optimal belief at time $t$. Then, for every reachable belief $\mu\in\mathcal R_t(\mf y^{[n]})$,
	\begin{equation}\label{eq:envelope-ineq-min}
		V_t(\mu)
		\;\le\;
		\big\langle U_t^{\mathrm R}(\cdot,\mf y^{[n]};\bar u),\,\mu\big\rangle,
	\end{equation}
	with equality at the optimal belief,
	\[
	V_t\big(\bar\mu_t^{\mf y^{[n]}}\big)
	=
	\big\langle U_t^{\mathrm R}(\cdot,\mf y^{[n]};\bar u),\,\bar\mu_t^{\mf y^{[n]}}\big\rangle.
	\]
	In particular, the map $\mu\mapsto\langle U_t^{\mathrm R}(\cdot,\mf y^{[n]};\bar u),\mu\rangle$ is an affine majorant of $V_t$ on the set of reachable beliefs, touching $V_t$ at $\bar\mu_t^{\mf y^{[n]}}$.
\end{lem}

\begin{rem}[Derivative along the submanifold of reachable beliefs]
	\label{rem:reachable-derivative}
	For fixed $(t,\mf y^{[n]})$, let $\mathcal R_t(\mf y^{[n]})$ denote the subset of $\mathcal P(\R^{d_x})$ consisting of all beliefs at time $t$ that are reachable under admissible partially observed controls started from $(t,\mf y^{[n]})$. We may view $\mathcal R_t(\mf y^{[n]})$ as a submanifold of $\mathcal P(\R^{d_x})$. Any one-parameter perturbation $(u^\theta)_\theta$ of the optimal control $\bar u$ generates a curve $\theta\mapsto\mu_t^\theta\in\mathcal R_t(\mf y^{[n]})$ with $\mu_t^0=\bar\mu_t^{\mf y^{[n]}}$. Assuming differentiability in $\theta$ at $0$, the signed measure
	\[
	\nu := \left.\frac{d}{d\theta}\mu_t^\theta\right|_{\theta=0}
	\]
	belongs to the tangent space $T_{\bar\mu_t^{\mf y^{[n]}}}\mathcal R_t(\mf y^{[n]})$ of this submanifold at $\bar\mu_t^{\mf y^{[n]}}$. Since each $\mu_t^\theta$ is a probability measure, $\nu$ has total mass zero,
	\[
	\int_{\R^{d_x}} \nu(dx) = 0.
	\]
	
	By Lemma~\ref{lem:envelope-ineq}, for all $\theta$ we have
	\[
	V_t(\mu_t^\theta)
	\;\le\;
	\big\langle U_t^{\mathrm R}(\cdot,\mf y^{[n]};\bar u),\,\mu_t^\theta\big\rangle,
	\qquad
	V_t\big(\bar\mu_t^{\mf y^{[n]}}\big)
	=
	\big\langle U_t^{\mathrm R}(\cdot,\mf y^{[n]};\bar u),\,\bar\mu_t^{\mf y^{[n]}}\big\rangle.
	\]
	Subtracting the equality at $\theta=0$ and dividing by $\theta>0$ yields
	\[
	\frac{V_t(\mu_t^\theta)-V_t(\bar\mu_t^{\mf y^{[n]}})}{\theta}
	\;\le\;
	\frac{\big\langle U_t^{\mathrm R}(\cdot,\mf y^{[n]};\bar u),\,\mu_t^\theta-\bar\mu_t^{\mf y^{[n]}}\big\rangle}{\theta}.
	\]
	Letting $\theta\downarrow 0$ and using differentiability of $\theta\mapsto V_t(\mu_t^\theta)$ and $\theta\mapsto\mu_t^\theta$ at $0$ gives the right-hand directional derivative
	\[
	D^+ V_t\big(\bar\mu_t^{\mf y^{[n]}};\nu\big)
	\;\le\;
	\int U_t^{\mathrm R}(x,\mf y^{[n]};\bar u)\,\nu(dx).
	\]
	Performing the same argument with $\theta<0$ yields the left-hand derivative inequality
	\[
	D^- V_t\big(\bar\mu_t^{\mf y^{[n]}};\nu\big)
	\;\ge\;
	\int U_t^{\mathrm R}(x,\mf y^{[n]};\bar u)\,\nu(dx).
	\]
	If $V_t$ admits a directional derivative along $\mathcal R_t(\mf y^{[n]})$ at $\bar\mu_t^{\mf y^{[n]}}$ in the direction $\nu$, then $D^+ V_t = D^- V_t$ and the two inequalities combine to give
	\[
	\left.\frac{d}{d\theta}V_t(\mu_t^\theta)\right|_{\theta=0}
	=
	\int U_t^{\mathrm R}(x,\mf y^{[n]};\bar u)\,\nu(dx),
	\qquad
	\nu \in T_{\bar\mu_t^{\mf y^{[n]}}}\mathcal R_t(\mf y^{[n]}).
	\]
	In particular, the kernel $x\mapsto U_t^{\mathrm R}(x,\mf y^{[n]};\bar u)$ represents the flat derivative of the restriction of $V_t$ to the submanifold $\mathcal R_t(\mf y^{[n]})$ at the optimal belief $\bar\mu_t^{\mf y^{[n]}}$. Since $\nu$ has zero total mass, adding any $x$-constant $C$ to $U_t^{\mathrm R}$ does not change the pairing since  $\langle C,\nu\rangle = C\int\nu(dx)=0$, thus the flat derivative is defined only up to an additive $x$-constant.
\end{rem}

\begin{propo}[Adjoint as flat derivative along the optimal belief path]
	\label{prop:adjoint-flat}
	Let $\mf y^{[n]}$ be a realized observation history, and let $(\bar u,\bar\mu,\bar U)$ be the globally optimal triple of Proposition~\ref{prop:aux-opt} for the partially observed problem started from $(t,\mf y^{[n]})$, so that~\eqref{eq:V-envelope-again} and the envelope inequality of Lemma~\ref{lem:envelope-ineq} hold. Assume moreover that $V_t\in\mathcal S^{1,1}(\mathcal P(\R^{d_x}))$ and is  differentiable at $\bar\mu_t^{\mf y^{[n]}}$. Then, for every $t\in(t_n,t_{n+1})$,
	\begin{equation}\label{eq:adjoint=flat}
		\frac{\delta V_t}{\delta\mu}\big(\bar\mu_t^{\mf y^{[n]}},x\big)
		= U_t^{\mathrm R}(x,\mf y^{[n]};\bar u),
	\end{equation}
	up to the usual additive constant in the flat derivative. Equivalently, we may write
	\[
	\frac{\delta V_t}{\delta\mu}\big(\bar\mu_t^{\mf y^{[n]}},x\big)
	= \bar U_t(x,\mf y^{[n]}).
	\]
\end{propo}

Equation~\eqref{eq:adjoint=flat} shows that, along the optimal belief
trajectory, the costate $\bar U_t(\cdot,\mf y^{[n]})$ coincides with the flat
derivative $\delta V_t/\delta\mu$ of the value functional in its measure
argument. In particular, $\bar U_t$ encodes the infinitesimal sensitivity of
the value to perturbations of the belief. This confirms that our Pontryagin
necessary conditions are consistent with the Hamilton-Jacobi-Bellman
characterisation of $V_t$ on the space of probability measures.

Thus, under the convexity and differentiability conditions, we obtain the infinite-dimensional envelope condition 
\[
\frac{\delta V_t}{\delta \mu}\Big(\bar\mu_{y^{[n]}_t};\,x\Big) \;=\; U_R\!\big(x,\,y^{[n]};\,\bar u\big)\,,
\] 
for all $x\in\R^{d_x}$ (up to an additive constant in $x$).  In other words, the adjoint $U_R(\cdot,\,y^{[n]};\,\bar u)$ coincides with the gradient of the value functional $V_t(\mu)$ at the optimal belief $\bar\mu_{y^{[n]}_t}$, just as $p_t = \nabla_x V(t,X^*_t)$ in the classical fully observed, deterministic case.

	\subsection{Separation principle}
	\label{subsec:separation}
	
	In linear–Gaussian models with quadratic costs, the separation principle says the optimal output–feedback controller is obtained by (i) estimating the state $\hat X_t$ from the available information (here $\mb F_t^Y$) and (ii) applying the full–state LQR law to $\hat X_t$ (certainty equivalence). This coincides with the known finite-horizon separation result (see, for example, \cite[Thm.~2.1]{Wonham1968}).
	
	\subsubsection{Assumptions (finite horizon, discrete observations)}
	We use \eqref{eq:POSOC_LQG_observation}–\eqref{eq:POSOC_LQG_cost_functional} with discrete observations. Separation holds on $[0,T]$ under:
	\begin{enumerate}
		\item $R\succ 0$ and $Q,Q_T\succeq 0$.
		\item \textbf{Linear–Gaussian with independent noises.} The state and observation models are linear–Gaussian; $W_t$, $\{S_n\}$, and $X_0$ are Gaussian and mutually independent (standard finite–horizon LQG assumptions\ \cite[Thm.~2.1]{Wonham1968}).
		\item \textbf{Exogenous observations (control–independent information).}
		\[
		\mathcal L\!\big(Y_{t_n}\mid X_{t_n},\mathbb F_{t_n^-}\big)
		= \mathcal N\!\big(CX_{t_n},\,R_{y,n}\big),\qquad R_{y,n}:=\beta_{t_n}\beta_{t_n}^\top,
		\]
		where the sampling times $\{t_n\}$ and covariances $\{R_{y,n}\}$ are fixed \emph{a priori} (do not depend on the control policy). Intuitively: choosing $\alpha$ cannot change future measurement quality. If this fails (information is \emph{endogenous}), the problem has a \emph{dual effect} and the optimal policy is \emph{not} certainty–equivalent; see \textbf{Bar–Shalom \& Tse (1974), Sec.~II} and the explicit networked–sensing example in \textbf{Ramesh–Sandberg–Bao–Johansson (2011), Sec.~III–B} \cite{BarShalomTse1974,RameshSandbergBaoJohansson2011}.
	\end{enumerate}
	
	\subsubsection{Optimal controller and estimator (decoupled)}
	With the above assumptions, the optimal control is
	\begin{equation}
		\label{eq:sep-optimal-law}
		\alpha_t^\star \;=\; -\,K(t)\,\hat X_t,
		\qquad K(t)=R^{-1}B^\top S(t),
	\end{equation}
	where $S:[0,T]\to\mb S_+^{d_x}$ solves the control Riccati ODE
	\begin{equation}
		\label{eq:control-riccati}
		-\,\dot S(t)=A^\top S(t)+S(t)A - S(t)B R^{-1}B^\top S(t) + Q,\qquad S(T)=Q_T,
	\end{equation}
	(cf.\ continuous–time LQR \cite{Kalman1960}). Let $\Sigma:=\sigma\sigma^\top$. Between observations $t\in(t_n,t_{n+1})$,
	\begin{equation}
		\label{eq:KF-ct-predict}
		\dl \hat X_t=\big(A\hat X_t+B\alpha_t\big)\dl t,\qquad
		\dot P_t = A P_t + P_t A^\top + \Sigma,
	\end{equation}
	and at $t_n$,
	\begin{equation}
		\label{eq:KF-discrete-update}
		\hat X_{t_n}^+ = \hat X_{t_n}^- + K_n\!\left(Y_{t_n}-C\hat X_{t_n}^-\right),\quad
		K_n = P_{t_n}^- C^\top \!\left(CP_{t_n}^-C^\top+R_{y,n}\right)^{-1},\quad
		P_{t_n}^+ = (I-K_n C)P_{t_n}^-,
	\end{equation}
	i.e., continuous–time prediction with discrete measurement updates \cite{KalmanBucy1961}. 
	The controller Riccati \eqref{eq:control-riccati} depends only on $(A,B,Q,R,Q_T)$, while the Kalman recursions \eqref{eq:KF-ct-predict}–\eqref{eq:KF-discrete-update} depend on $(A,C,\Sigma)$ and the fixed $\{R_{y,n}\}$.
	
	\subsubsection{When exogeneity fails: ”dual effect and consequent lossof separation}
If the observation law can be influenced by control (for example, if $\{\beta_{t_n}\}$ is a decision variable), then information is \emph{endogenous}: future information quality depends on the control policy. This induces the \emph{dual effect}, and, in general, the optimal controller is not certainty–equivalent: estimation and control cannot be designed independently (the separation principle fails). See \cite[Sec.~II]{BarShalomTse1974} for the precise “certainty equivalence $\Leftrightarrow$ no dual effect’’ statement and \cite[Sec.~III–B]{RameshSandbergBaoJohansson2011} for a concrete example where the estimator covariance depends on control.

\section{Numerical approach}
\label{sec:numerical_approach}

We present a particle–based scheme to approximate the pathwise forward–backward system \eqref{eq:aux_system}. 
The method searches for a fixed point in a parametric representation of the needle–induced pathwise value function $\overline U$. 
Given a candidate $\hat p(\cdot;\theta)$, we (i) extract policies by minimizing belief–weighted Hamiltonians, (ii) propagate particles and update beliefs forward in time, and (iii) regress from simulated rollouts to update $\theta$. 
At convergence the iteration satisfies the fixed–point relation $\theta^\star=\mathcal T(\theta^\star)$, yielding a self–consistent locally optimal control law.

\subsection{Particle representation of beliefs}
\label{subsec:particles}

Let $\{X_{0}^{(m)}\}_{m=1}^M\sim\mu_0$ be $M$ particles i.i.d sampled from the initial distribution. 
The empirical measure $\mu_t^M:=\tfrac{1}{M}\sum_{m=1}^M \delta_{X_t^{(m)}}$ approximates the evolving belief. 
Between observation times, particles follow Euler–Maruyama:
\[
X^{(m)}_{t+\Delta t}
= X^{(m)}_t + b\!\big(t,X^{(m)}_t,\alpha_t\big)\,\Delta t
+ \sigma\!\big(t,X^{(m)}_t,\alpha_t\big)\,\Delta W_t^{(m)}.
\]
At observation times $t_n$, particles are reweighted by the likelihood $\pi_n$ and resampled, when needed, to approximate the Bayesian posterior.

\subsection{Parametrization of the pathwise value and a finite observation window}
\label{subsec:parametrization}

On each inter–observation slab $[t_n,t_{n+1})$, we approximate the auxiliary-cost-to-go  as
\[
\overline U_t^{[n]}(x,\mf y^{[n]})\;\approx\;\hat p_t^{\theta}\!\big(x,\;\mf y^{[n]}\big),
\]
with $\theta$ the parameters of the chosen family (polynomials, kernels, or neural networks). 
A direct parametrization in the full observation vector $\mf{Y}^{[n]}=(Y_{t_0},\dots,Y_{t_n})$ suffers from a steadily growing input dimension as $n$ increases, which degrades sample efficiency and complicates regression.
To control this growth, we \emph{by design} condition only on the most recent $K$ observations. 
Define the sliding window
\[
\mf Z^{[n]} \;:=\; \big(Y_{t_{n-K+1}},\,\dots,\,Y_{t_n}\big)\in\R^{K d_y}
\quad\text{(with the obvious truncation if $n<K-1$),}
\]
and rewrite the approximation as
\[
\overline U_t^{[n]}(x,\mf z^{[n]})\;\approx\;\hat p_t^{\theta}\!\big(x,\;\mf z^{[n]}\big).
\]
For notational brevity we denote the window simply by $\mf z$ when the time index is clear from context. 
The hyperparameter $K$ trades off statistical efficiency (larger $K$ captures longer memory) against computational tractability and variance in regression.

\subsection{Policy extraction via regression}
\label{subsec:policy}

Policies are extracted from $\hat p^\theta$ by minimizing conditional expectations given the window $\mf z$.
For the continuous control,
\[
\bar{\alpha}_{t}(\mf z)
\;\in\;\argmin_{\alpha\in\mcal A}\;
\mb E\!\Big[f_t(X_t,\alpha)
+\mcal G_\alpha \hat p_t^\theta\!\big(X_t,\mf z\big)
\;\big|\;\mf z^{[n]}=\mf z\Big],
\]
and for the discrete control at observation times,
\[
\bar{\beta}_{t_n}(\mf z^-)
\;\in\;\argmin_{\beta\in\mcal B}\;
\mb E\!\Big[c_{t_n}\!\big(X_{t_n^-},\beta\big)
+\hat p_{t_n}^\theta\!\big(X_{t_n^-},\mf z\big)
\;\big|\;\mf z^{[n-1]}=\mf z^-\Big],
\]
where $\mf z^-:=\mf z^{[n-1]}$ denotes the pre–observation window at $t_n$. 
Both conditional expectations are estimated by $L^2$ regression with particles using simulated pairs $\{(X_t^{(m)},\mf z^{(m)})\}_{m=1}^M$ resulting in ($\alpha_{t}^{(m)}, \beta_{t_n}^{(m)}$). 
In special cases (e.g., LQG systems) the minimizers admit closed forms.

\subsection{Parameter update by regression}
\label{subsec:regression}

After simulating trajectories under the current parameters, define the pathwise costs
\[
P_t^{(m)} := \int_t^T f_\tau\!\big(X_\tau^{(m)},\alpha_\tau^{(m)}\big)\,\dl \tau
+ \sum_{i:\,t_i\ge t} c_{t_i}\!\big(X_{t_i^-}^{(m)},\beta_{t_i}^{(m)}\big)
+ g\!\big(X_T^{(m)}\big).
\]
Update $\theta$ by least squares:
\[
\theta^{\ell+1} \;\in\; \argmin_\theta \sum_{m=1}^M
\Big(\hat p_t^{\theta}\!\big(X_t^{(m)},\mf z^{(m)}\big) - P_t^{(m)}\Big)^2,
\]
where $\ell$ is the outer iteration index. The parametrization $\theta$  depends on time $t$. For notational simplicity, we suppress this dependence and write $\theta$ instead of $\theta(t)$. In particular, under a time-grid discretization $\{t_n\}_{n=0}^{N+1}$, we allow for a distinct parameter vector $\theta_n := \theta(t_n)$ at each grid point.

\label{subsec:algorithm}
\RestyleAlgo{ruled}
\SetKwComment{Comment}{/* }{ */}
\begin{algorithm}[h!]
	\caption{a article fixed-point method}
	\KwIn{$\mu_0, \rho_{0}, \theta_0, t_0, T, M, N, N_{o}, n_{\text{iter}}, d_x, d_y, K$}
	\tcc{Initialization}
	Sample $\{X_{t_0}^{(m)}\}_{m=1}^M\sim \mu_0$\;
	Sample $\{Y_{t_0}^{(m)}\}_{m=1}^M\sim \rho_0$ \Comment*[r]{if available}
	$\ell \gets 0$\;
	\While{stopping criterion not met}{
		\For{$n=0$ \KwTo $N_{o}$}{
			\tcc{Observation and control}
			Sample $\{Y_{t_n}^{(m)}\}_{m=1}^M$ using noise level parameter $\bar\beta_n$\;
			Form windows $\{\mf z^{[n],(m)}\}$ of length $K$ from $\{Y_{t_i}^{(m)}\}_{i\le n}$\;
			Estimate $\bar\alpha_t(\mf z)$ by $L^2$ regression of the conditional expectation\;
			Evolve particles forward under $\bar\alpha_t$ on $[t_n,t_{n+1})$\;
			Estimate $\bar\beta_{n+1}(\mf z)$ by $L^2$ regression at $t_{n+1}$\;
		}
		\tcc{Regression update}
		Compute pathwise costs $\{P_t^{(m)}\}_{m=1}^M$ and set $P_T^{(m)}=g(X_T^{(m)})$\;
		Update $\theta^{\ell+1}\gets \argmin_\theta \sum_{m=1}^M
		\big(\hat p^\theta_t(X_t^{(m)},\mf z^{(m)})-P_t^{(m)}\big)^2$\;
		$\ell\gets \ell+1$\;
	}
\end{algorithm}

\subsection{Numerical example: LQG under partial observations}
\label{subsec:lqg_numerical_example}

We demonstrate the particle fixed–point scheme on the partially observed LQG model \eqref{eq:POSOC_LQG_observation}.
All experiments use the numerical pipeline of Section~\ref{sec:numerical_approach} (Euler–Maruyama propagation,
Bayesian reweighting at observation times, $L^2$ regression for conditional expectations, and the finite
observation window $\mf z$ for tractable regression).

We consider two cases:
\begin{itemize}
	\item \textbf{Low-dimensional (1D)} for transparent comparisons with closed-form LQG/separation solutions.
	\item \textbf{High-dimensional} ($d_x=10$) where grid-based HJB methods are infeasible.
\end{itemize}
For each scale we evaluate two observation-noise regimes:
\begin{enumerate}[ (A)]
	\item \textbf{Fixed observation noise} ($\beta_{t_n}\equiv \varepsilon$): 
	\(
	Y_{t_n}=C X_{t_n}+\varepsilon\,\xi_n,\quad \xi_n\sim\mathcal N(0,I_{d_y}),
	\)
	with prescribed $\varepsilon>0$.
	\item \textbf{Control-dependent observation noise} (optimized $\beta$): a discrete control modulates the
	noise level and incurs an observation cost \eqref{eq:POSOC_LQG_cost_functional}.
\end{enumerate}

\subsubsection{Low-dimensional (1D) results}

\paragraph{Regime A (fixed observation noise).}

In the LQG example we vary $N_o$ while keeping all other hyperparameters fixed and we also compare to the FOSOC (Fully Observed Stochastic Optimal Control).
\begin{table}[H]
	\centering
	\caption{ \textbf{ LQG } Expected cost [with 95\% confidence interval] for particle method vs separation principle cost for different numbers of observations $N_o$, $d_x{=}1$. Observation memory fixed to $K{=}1$. Parameters : $A{=}{-}0.25$, $B{=}C{=}1.0$, $\sigma{=}0.5$, $ Q {=} R{=}Q_T{=}2.0$, $\ve {=} 0.1$, $T{=}1.0$, $M_{\text{eval}}{=}10^5$, $M_{\text{train}}{=}500$.}
	\label{tab:1d_avg_costs_alt}
	\begin{tabular}{c c c}
		\toprule
		$N_o$ & Particle ([95\% CI]) & Separation \\
		\midrule
		1   & $1.37\,[1.35,\,1.39]$ & $1.373$ \\
		5   & $1.15\,[1.13,\,1.16]$ & $1.150$ \\
		10  & $1.09\,[1.08,\,1.11]$ & $1.095$ \\
		30  & $1.06\,[1.04,\,1.07]$ & $1.051$ \\
		\midrule
		$\text{FOSOC}$ & \multicolumn{2}{c}{$1.024$} \\
		\bottomrule
	\end{tabular}
\end{table}

Table~\ref{tab:1d_avg_costs_alt} shows that our particle-based method achieves an expected cost nearly identical to the optimal cost obtained by the  separation principle across varying numbers of observations $N_o$, with discrepancies lying within the statistical confidence intervals. Increasing $N_o$ monotonically reduces the expected cumulative cost and narrows the gap with the FOSOC reference. This confirms that, in the exogenous–observation setting where the separation principle holds, the algorithm effectively recovers the optimal strategy.

However, the separation–principle curve is the partial-information benchmark in the LQG setup, it conditions on the \emph{entire} observation history, whereas our particle fixed-point scheme uses a finite window with memory $K=1$, i.e., policies and value depend only on the most recent observation. Consequently, discrepancies between the two are not expected to vanish strictly as $N_o$ increases and need not be monotone in $N_o$. The observed gaps can be explained by (i) Monte Carlo error (finite $M_{\text{train}}$, $M_{\text{eval}}$), (ii) time-discretization bias in the forward–backward iteration (finite $\Delta t$), and (iii) statistical/approximation error from the $L^2$ regression used for policy extraction, \emph{plus} (iv) deliberate information truncation from using $K=1$ instead of the full observation history. The case $N_o=30$ is particularly informative: the separation solution has assimilated many past measurements, whereas our controller discards all but the latest one. Any residual gap there should thus be attributed at least in part to the windowing design rather than to a failure of the scheme. Within these constraints, the table indicates that the algorithm converges to the separation benchmark up to a small, well-understood tolerance.

\begin{figure}[h!]
	\centering
	\includegraphics[width=\textwidth]{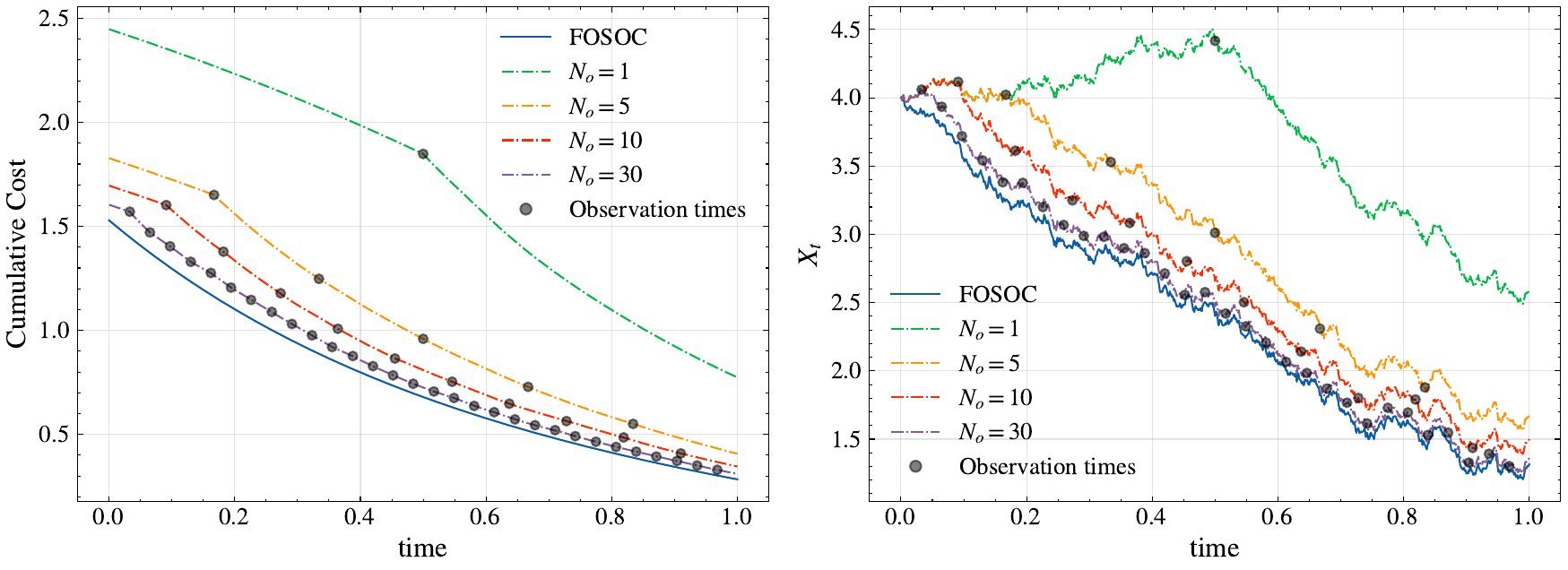}\hfill
	\caption{\textbf{ LQG fixed noise $\varepsilon$, $d_x{=}1$ .}
		Left: expected cost-to-go vs.\ time for different $N_o$ and the FOSOC benchmark.
		Right: example state paths; tighter regulation with larger $N_o$. Parameters : $A{=}0.25$, $B{=}C{=}1.0$, $\sigma{=}0.5$, $ Q {=} R{=}Q_T{=}2.0$, $\ve {=} 0.1$, $T{=}1.0$, $M_{\text{eval}}{=}10^5$, $M_{\text{train}}{=}2000$}
	\label{fig:1d_fixed_noise_multiple_N_o}
\end{figure}

In Fig.~\ref{fig:1d_fixed_noise_multiple_N_o} (left), the cost-to-go curves tighten and approach the fully observed benchmark as the observation frequency ($N_o$) increases, consistent with Table~\ref{tab:1d_avg_costs_alt}. The right panel shows sample trajectories for varying $N_o$ with initial state $x_0=4.0$ and drift $A=0.25$. Uncontrolled, the state drifts away from zero, incurring higher cost; under the optimal controller it is driven toward zero. To isolate the effect of observation timing, all runs share identical Wiener increments. Earlier observations improve the state estimate and enable earlier, stronger corrections, bringing the trajectory closer to zero sooner.

\paragraph{Regime B (noise as a control).}
We jointly optimize $\alpha$ and $\beta$ with observation cost :  $$\sum_{n=1}^{N_o} c_n(\beta_{t_n}) = \sum_{n=1}^{N_o} \sum_{i = 1}^{d_y} \frac{\kappa_{n,i}}{\beta_{t_n,i}} =  \sum_{n=1}^{N_o} \trace(\diag(\kappa_n) \diag^{-1}(\beta_{t_n})) $$
exhibiting the classic exploration–information tradeoff:
smaller $\beta$ improves information but increases observation cost.
\begin{figure}[H]
	\centering
		\includegraphics[width=0.65\textwidth]{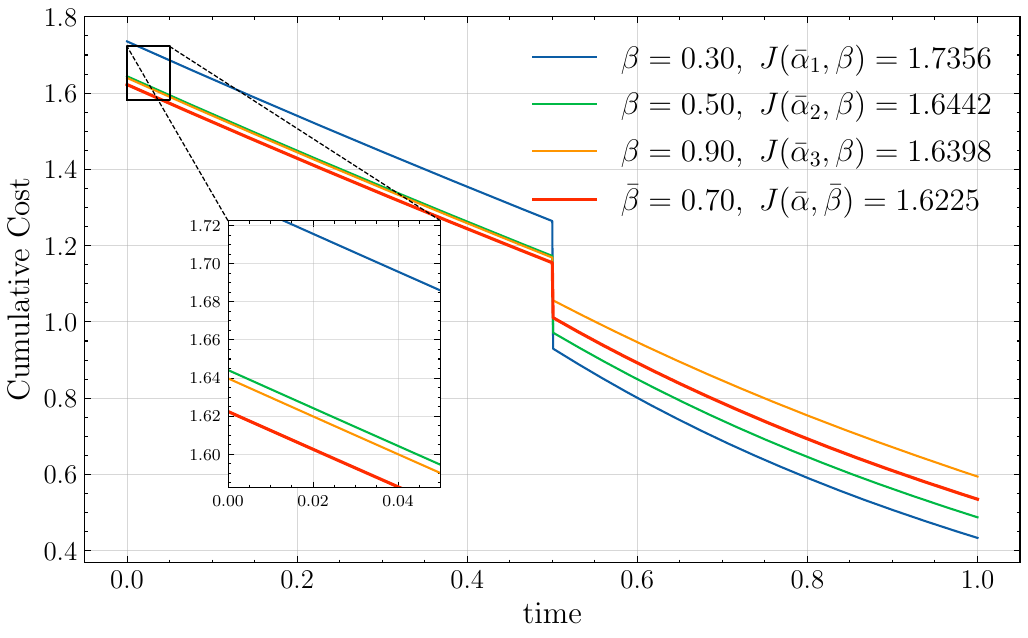}
	\caption{\textbf{ LQG controlled noise $\bar\beta$, $d_x=1$, $N_o = 1$.}
		Expected cost-to-go for different non-optimal noise levels ($\beta = 0.3,0.5,0.9$),  vs. optimal noise level $\bar{\beta}$ (scalar). Observation memory fixed to $K{=}1$. Parameters : $A{=}{-}0.25$, $B{=}C{=}1.0$, $\sigma{=}0.5$, $ Q {=} R{=}Q_T{=}2.0$, $\ve {=} 0.1$, $\kappa_{1}{=}0.1$, $T{=}1.0$, $M_{\text{eval}}{=}10^5$, $M_{\text{train}}{=}10^4$. }
	\label{fig:1d_controlled_noise_1obs}
\end{figure}
The change in ordering of the curves after the observation time $t_1=0.5$
has a clear interpretation: once the discrete sensing cost has been paid,
we are simply running an optimal LQG controller for the problem for $t>t_1$,
starting from the posterior produced by the chosen $\beta$.  Different
$\beta$’s induce different posteriors (information quality), and hence
different optimal costs, which is exactly what the dynamic
programming principle predicts.

We also showcase results with a higher number of observation, where the discrete control $\beta$ becomes a function of the previous observations :  
\begin{figure}[H]
	\centering
	\includegraphics[width=0.55\textwidth]{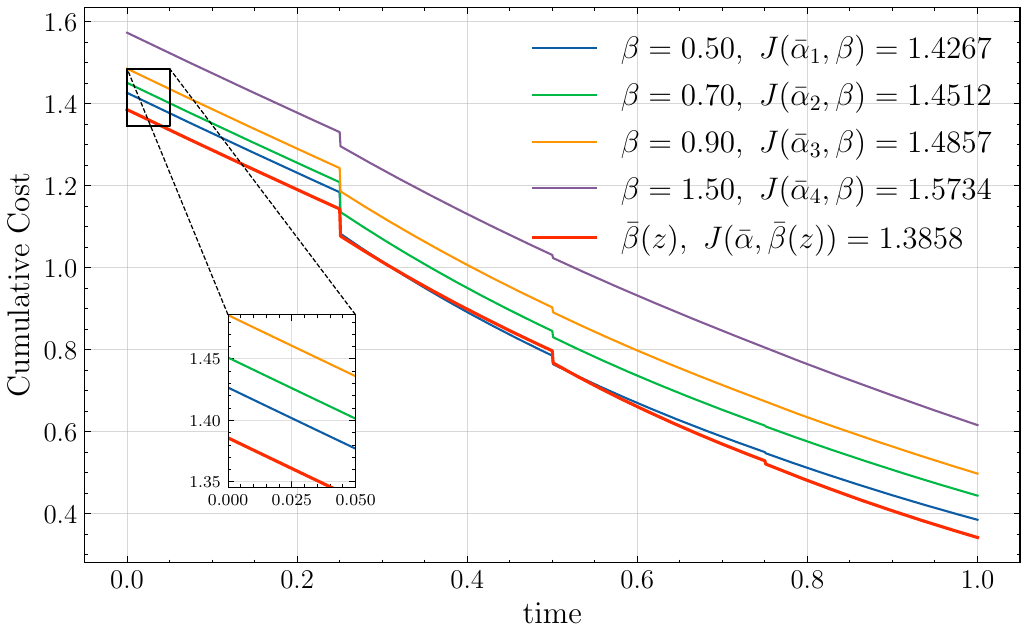}
	\includegraphics[width=0.35\textwidth]{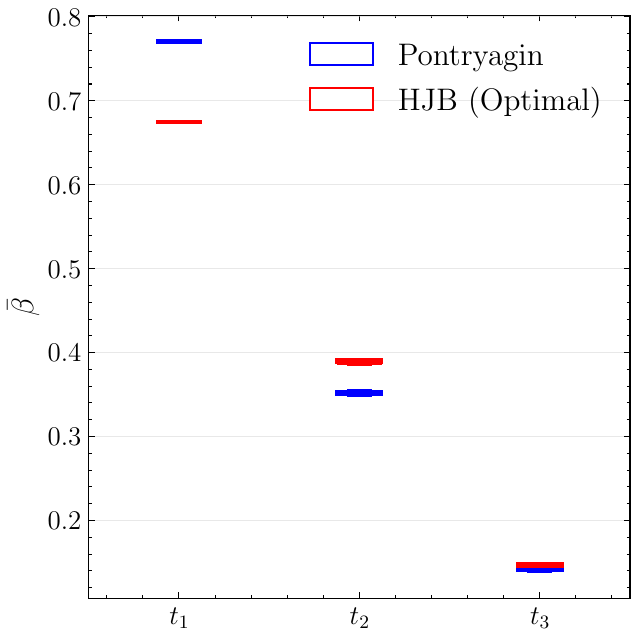}
	\caption{\textbf{ LQG  controlled noise $\bar\beta$, $d_x=1$, $N_o = 3$.}
	Left : Expected cost-to-go for different non-optimal noise levels ($\beta = 0.5,0.7,0.9, 1.5$), vs. optimal noise level function $\bar{\beta}(z)$.
	Right :  Optimal noise levels $\bar{\beta}$ at each observation step for realizations of the process.
	Observation memory fixed to $K{=}2$. Parameters : $A{=}{-}0.25$, $B{=}C{=}1.0$, $\sigma{=}0.5$, $ Q {=} R{=}Q_T{=}2.0$, $\kappa{=}[0.05,0.01,0.001]$, $T{=}1.0$, $M_{\text{eval}}{=}10^5$, $M_{\text{train}}{=}10^4$.}
	\label{fig:1d_controlled_noise_3obs}
\end{figure}

Figure \ref{fig:1d_controlled_noise_3obs} illustrates the effect of controlling the noise level when three observations are available.
On the left, we plot the cumulative cost for several fixed noise levels $\beta \in \{0.5,0.7,0.9, 1.5\}$ and for the control $\bar\beta(z)$ learned from our Pontryagin–type formulation.
For all times, the adaptive control $\bar\beta(z)$ yields a strictly smaller expected cost-to-go than any of the constant choices, leading to the lowest cost.

The right panel compares, at each observation time $t_1,t_2,t_3$, the optimal discrete noise levels obtained with our method (blue, “Pontryagin”) to the benchmark values computed previously from the HJB equation (red) as in \cite{bayer2024continuoustimestochasticoptimal}.  The two sets of controls almost coincide, showing that our parametric ansatz for the adjoint and for $\beta$ is sufficiently expressive to reproduce the HJB-based solution.  The  discrepancies that remain can be attributed to the  ansatz choice  and to the finite number of particles used in the forward–backward simulations. This trend is also consistent with our cost structure: since the coefficients $\kappa = [0.05,0.01,0.001]$ decrease over observation times, later observations are cheaper to use, so one naturally expects smaller optimal noise levels $\beta_{t_n}$ when the associated cost is smaller.

\subsubsection{High-dimensional results}
We increase the state dimension (to $d_x=10$) while retaining linear–Gaussian dynamics
and quadratic costs. Grid-based HJB methods are omitted due to infeasibility; separation-principle
baselines remain available in Regime~A.
\begin{figure}[H]
	\centering
	\includegraphics[width=0.55\textwidth]{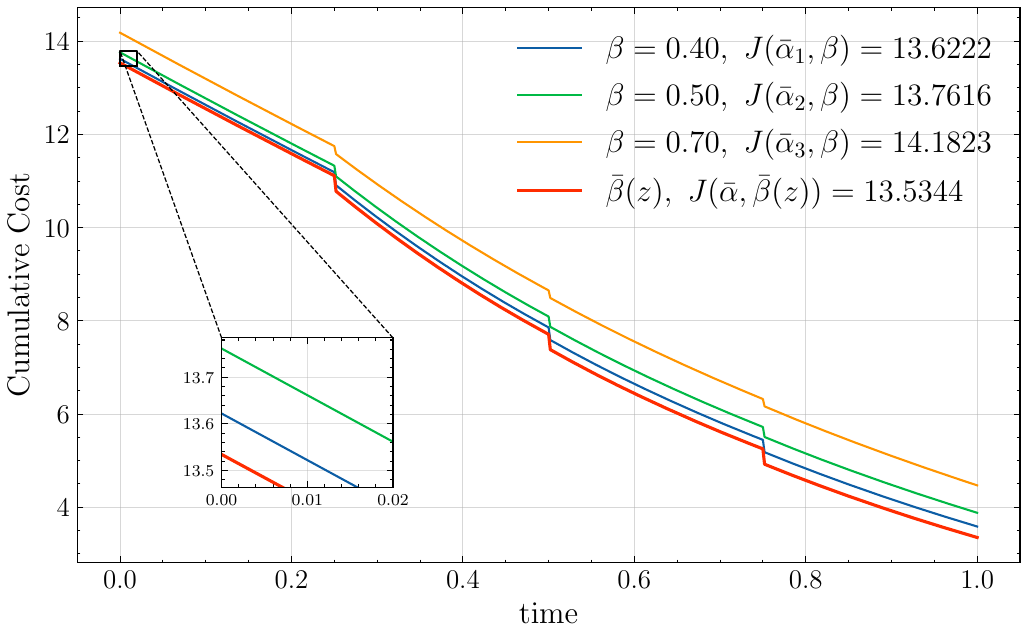}
	\caption{\textbf{ LQG controlled noise $\bar\beta$, $d_x=10$, $N_o = 3$.}
		Expected cost-to-go for different non-optimal noise levels ($\beta = 0.4,0.5,0.7$),  vs. optimal noise level function $\bar{\beta}(z)$. Observation memory fixed to $K{=}1$. Parameters : $A{=}{-}0.25 I_{d_x}$, $B{=}C{=}1.0 I_{d_x} $, $\sigma{=}0.5 I_{d_x}$, $ Q {=} R{=}Q_T{=}2.0 I_{d_x}$, $\kappa{=}[0.1,0.1,0.1]$, $T{=}1.0$, $M_{\text{eval}}{=}10^5$, $M_{\text{train}}{=}10^4$.}
	\label{fig:multid_controlled_noise_3obs}
\end{figure}
Figure \ref{fig:multid_controlled_noise_3obs} illustrates the high-dimensional ($d_x = 10$) performance of our Pontryagin-based scheme with three observation opportunities at $t_1$, $t_2$ and $t_3$. The adaptive noise policy $\bar{\beta}(z)$ (solid red curve) attains a consistently lower expected cost-to-go than any fixed noise level, confirming that the algorithm effectively balances exploration and exploitation even in higher dimensions. Each jump in the cost trajectory at $t_1$, $t_2$ and $t_3$ reflects the instantaneous observation cost and the subsequent update in expected future cost due to new information. The adaptive controller achieves a cost of $J = 13.5340$, outperforming the best fixed noise level ($\beta = 0.40$, $J = 13.6222$). This demonstrates the scalability of our method to moderate state dimensions where traditional grid-based HJB solvers are intractable.
\subsection{Numerical Example: Non-LQG under Partial Observations}
\label{subsec:non_lqg_numerical_example}
We now introduce a benchmark problem to illustrate the behavior of the proposed framework
in a setting that lies outside the LQG class due to non-quadratic running costs that was introduced in \cite{tottori_forward-backward_2023} that we adapted to discrete observation setup.
The dynamics of the controlled process $X_t$ are given by
\begin{equation}
	\label{eq:nonlqg_dynamics}
	\begin{aligned}
		\dl X_t &= \alpha_t \, \dl t \,+\, \sigma \, \dl W_t, \\
		X_0 &\sim \mcal{N}(0,\,0.01),
	\end{aligned}
\end{equation}
where $X_t \in \mb{R}^{d_x}$, $\alpha_t \in \mb{R}^{d_\alpha}$ is the control,
$\sigma \in \mb{R}^{d_x}$ is the diffusion parameter, and $W_t$ is a standard Wiener process. 
We assume access only to discrete, noisy measurements of the state process. 
For a prescribed number of observation dates $N_o$ on $[0,T]$, the observations are given by
\begin{equation}
	\label{eq:nonlqg_observations}
	Y_{t_n} = C \, X_{t_n} + \varepsilon \, \xi_n, 
	\qquad \xi_n \sim \mcal{N}(0,I_{d_y}),
	\quad n=1,\dots,N_o,
\end{equation}
where $C \in \mb{R}^{d_y \times d_x}$ is the observation matrix and $\varepsilon > 0$ controls the noise level.The goal is to find the optimal control policy that minimizes
\begin{equation}
	\label{eq:nonlqg_cost}
	\begin{aligned}
		J(\alpha) &= \mb{E}\Bigg[ \int_0^T \Big( P(t,X_t) \,+\tfrac{1}{2} \,X_t^{\top} \, Q \, X_t + \tfrac{1}{2} \alpha_t^{\top} R \alpha_t \Big) \, \dl t 
		\,+\, \tfrac{1}{2}(X_T - x_\star)^{\top} Q_T (X_T - x_\star) \Bigg],
	\end{aligned}
\end{equation}
where $R \in \mb{R}^{d_\alpha \times d_\alpha}$, $Q, Q_T \in \mb{R}^{d_x \times d_x}$ are positive semi-definite,
and $P(t,x)$ encodes obstacle penalties that render the problem non-quadratic and $x_\star \in \mb{R}^{d_x}$ is the desired target state at the terminal time  .
We consider two scales:
\begin{itemize}
	\item \textbf{Low-dimensional (1D)} for transparent comparisons with closed-form LQG/separation solutions.
	\item \textbf{High-dimensional} ($d_x=10$) where grid-based HJB methods are infeasible.
\end{itemize}
For each scale we evaluate the following observation-noise regime:
\begin{enumerate}
	\item \textbf{Fixed observation noise} ($\beta_{t_n}\equiv \varepsilon$): 
	\(
	Y_{t_n}=C X_{t_n}+\varepsilon\,\xi_n,\quad \xi_n\sim\mathcal N(0,I_{d_y}),
	\)
	with prescribed $\varepsilon>0$.
\end{enumerate}
\paragraph{One-dimensional example.}
For $d_x=1$, the obstacle region is defined by
\begin{equation}
	P(t,x) =
	\begin{cases}
		1000, & 0.3 \leq t \leq 0.6 \;\; \text{and} \;\; 0.1 \leq |x| \leq 2, \\
		0,    & \text{otherwise}.
	\end{cases}
\end{equation}
The system want therefore to avoid the band $0.1 \leq |x| \leq 2$ during the time window $[0.3,0.6]$.
\begin{figure}[h!]
	\centering 
	\includegraphics[width=0.8\textwidth]{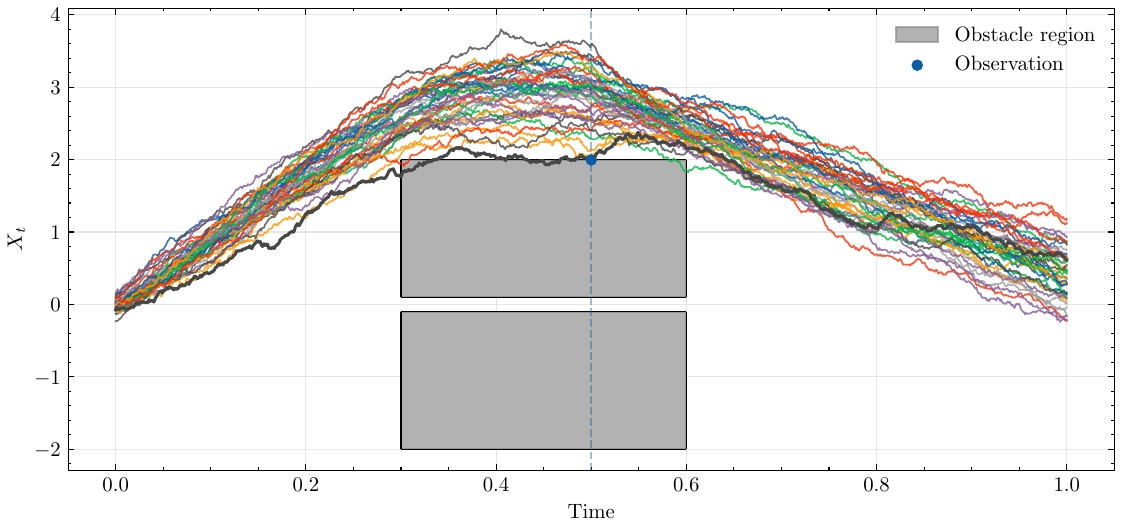}
	\caption{\textbf{NLQG  controlled trajectories, $d_x=1$,  $N_o = 1$.} Controlled sample trajectories of the state $X_t$. Observation memory fixed to $K{=}1$. Parameters : $C{=}1.0$, $\sigma{=}0.5$, $ Q {=} 0 $, $R{=}2.0$, $Q_T{=}20.0$, $x_{\star} = 0$, $\ve {=} 0.1$, $T{=}1.0$, $M_{\text{train}}{=}500$.}
	\label{fig:1d_nlqg_traj_highlighted}
\end{figure}
Figure~\ref{fig:1d_nlqg_traj_highlighted} displays sample trajectories of the state $X_t$ for the 1D non-LQG problem with a single observation ($N_o = 1$) and a non-convex obstacle penalty. The obstacle region (shaded in grey) is active during the time window $[0.3, 0.6]$ and the state band $0.1 \leq |x| \leq 2$. The trajectories, generated by our Pontryagin-based policy, successfully avoid the high-penalty region, with only minor incursions that are quickly corrected. Despite the non-quadratic cost structure and partial information, the algorithm produces a policy that effectively regulates the state toward the target $x_\star = 0$ by the terminal time, demonstrating its ability to handle complex, non-convex constraints under partial observations.
The trajectory highlighted in black demonstrates the informational value of the observation. Prior to the observation time, the controller operates with limited certainty about the state, and the trajectory drifts slightly inside the obstacle region. At the moment of observation (in blue), a noisy measurement is obtained, which allows the belief to be updated via Bayes’ rule. Immediately afterward, the controller leverages this new information to steer the trajectory decisively away from the high-penalty zone. This visible change in the trajectory after the observation illustrates how partial-observation control actively uses sparse measurements to correct course and avoid costly regions.

For $d_x=2$, the state $X_t = (X_t^{(1)}, X_t^{(2)})^\top$ evolves under
\eqref{eq:nonlqg_dynamics}--\eqref{eq:nonlqg_cost}.
We penalize trajectories that pass through a radial annulus around the origin during a fixed time window.
Geometrically, in $(t,x_1,x_2)$-space this forms a hollow cylinder (“tube”) active only for $t\in[0.3,0.6]$.
Let $r_{\rm in}<r_{\rm out}$ be the inner/outer radii and $x=(x_1,x_2)$. The running penalty is

$$
P(t,x)=
\begin{cases}
	1000, & 0.3\le t\le 0.6,\quad r_{\rm in}\le \|x\|_2\le r_{\rm out},\\[3pt]
	0, & \text{otherwise}.
\end{cases}
$$

Equivalently,

$$
P(t,x)=1000\;\mathbf{1}_{[0.3,0.6]}(t)\;\mathbf{1}_{[r_{\rm in},\,r_{\rm out}]}\!\big(\|x\|_2\big).
$$

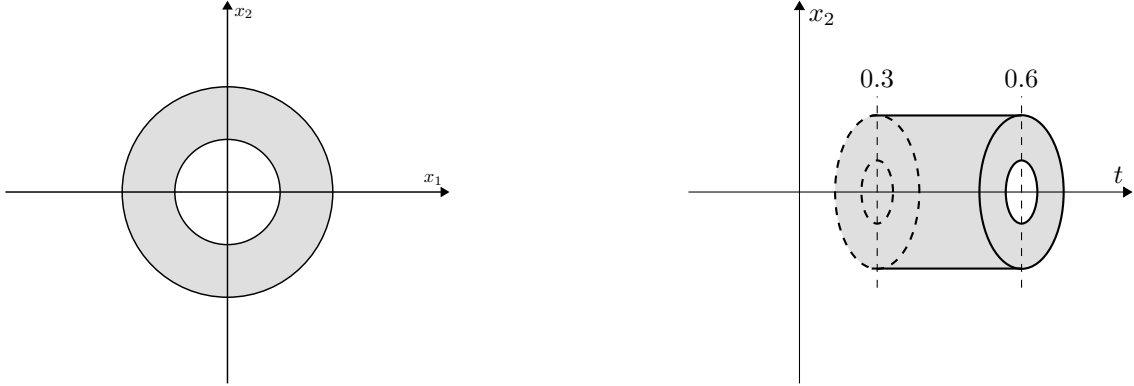
\begin{figure}[h]
	\centering
	\begin{minipage}{0.45\textwidth}
		\centering
		\begin{tikzpicture}[scale=0.70]
			\begin{axis}[
				myaxis,
				xmin=-3, xmax=3,
				ymin=-3, ymax=3,
				xlabel={$x_1$},
				ylabel={$x_2$},
				]
				\path[obsfill, even odd rule]
				(axis cs:0,0) circle[radius=2cm]
				(axis cs:0,0) circle[radius=1cm];
				
				\draw[obsline] (axis cs:0,0) circle[radius=2cm];
				\draw[obsline] (axis cs:0,0) circle[radius=1cm];
			\end{axis}
		\end{tikzpicture}
	\end{minipage}
	\hfill
	\begin{minipage}{0.45\textwidth}
		\centering
		\begin{tikzpicture}[scale=0.70]
			\begin{axis}[
				myaxis,
				xmin=-1, xmax=3,
				ymin=-5, ymax=5,
				xlabel={$t$},
				ylabel={$x_2$},
				]
				
				\addplot[name path=top, draw=none, samples=2, domain=0.7:2]{ 2};
				\addplot[name path=bot, draw=none, samples=2, domain=0.7:2]{-2};
				\addplot[obsfill, draw=none]
				fill between[of=top and bot];
				
				\path[obsfill, draw=none]
				(axis cs:0.7,0) ellipse [x radius=0.8cm, y radius=1.46cm];
				
				\path[obsfill, draw=none, even odd rule]
				(axis cs:2,0) ellipse [x radius=0.8cm, y radius=1.46cm]
				(axis cs:2,0) ellipse [x radius=0.3cm, y radius=0.6cm];
				
				\path[fill=white, draw=none]
				(axis cs:2,0) ellipse [x radius=0.3cm, y radius=0.6cm];
				
				\addplot[obsline, samples=2, domain=0.65:2]{ 2};
				\addplot[obsline, samples=2, domain=0.65:2]{-2};
				
				\draw[obsline] (axis cs:2,0) ellipse [x radius=0.8cm, y radius=1.46cm];
				\draw[obsline] (axis cs:2,0) ellipse [x radius=0.3cm, y radius=0.6cm];
				
				\draw[hideline] (axis cs:0.7,0) ellipse [x radius=0.8cm, y radius=1.46cm];
				\draw[hideline] (axis cs:0.7,0) ellipse [x radius=0.3cm, y radius=0.6cm];
				\draw[dashed] (0.7,-2.5) -- (0.7,2.5) node[above] {$0.3$};
				\draw[dashed] (2.0,-2.5) -- (2.0,2.5) node[above] {$0.6$};
			\end{axis}
		\end{tikzpicture}
	\end{minipage}
	\caption{Two-dimensional obstacle region (shaded) active for \(t\in[0.3,0.6]\).}
	\label{fig:nonlqg_obstacles_2d}
\end{figure}

\emph{Extension to $d_x=n$.}
Replace $x\in\mathbb{R}^2$ by $x\in\mathbb{R}^{d_x}$ and keep the same form with the Euclidean norm:

$$
P(t,x)=1000\;\mathbf{1}_{[t_{\min},t_{\max}]}(t)\;\mathbf{1}_{[r_{\rm in},\,r_{\rm out}]}\!\big(\|x\|_2\big),\qquad x\in\mathbb{R}^{d_x}.
$$
which yields a time-gated spherical shell in $\mathbb{R}^{d_x}$. 
\begin{figure}[h!]
	\centering
	\includegraphics[width=0.8\textwidth]{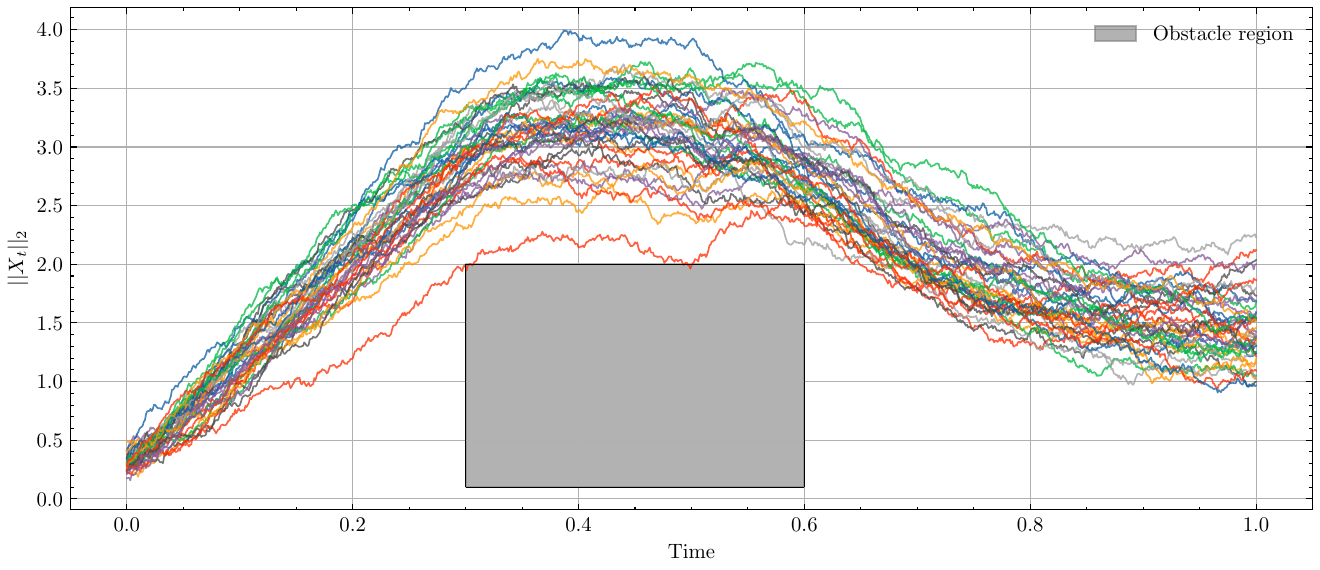}
	\caption{\textbf{$d_x=10$, NLQG  controlled trajectories, $N_o = 1$.} Controlled sample trajectories of the state $X_t$. Observation memory fixed to $K{=}1$. $C{=}1.0 I_{d_x}$, $\sigma{=}0.5 I_{d_x} $, $ Q {=} 0 $, $R{=}2.0 I_{d_x}$, $Q_T{=}20.0 I_{d_x} $, $x_{\star} = 0$, $\ve {=} 0.1$, $T{=}1.0$, $M_{\text{train}}{=}1200$.}
	\label{fig:10d_nlqg_traj}
\end{figure}

Figure \ref{fig:10d_nlqg_traj} illustrates the scalability of our approach to a high-dimensional non‑LQG setting ($d_x = 10$) with a spherical‑shell obstacle. The plot shows the Euclidean norm \(\|X_t\|_2\) of several sample trajectories under the learned Pontryagin‑based policy. The obstacle region (shaded in gray) is active for $t \in [0.3, 0.6]$ and for norms between an inner and outer radius. Despite the state’s high dimensionality and the non‑convex geometry of the penalty, the controller successfully avoid the penalty region throughout the critical time window. The single observation  occurs precisely at $t_1 = 0.5$, immediately afterward, the trajectories exhibit a coordinated descent that avoids the high‑cost region. This demonstrates that the algorithm extracts actionable information from a single noisy measurement and translates it into effective avoidance behavior even in a space where traditional grid‑based HJB methods are intractable.

In all experiments in this section, the outer fixed-point loop for the
value-function parameters was iterated until the improvement in the
estimated cost fell below a prescribed tolerance. Empirically, the
parameter iterates converged to a stable value $\theta^\star$. Moreover,
when we repeated the procedure from several distinct initial
parametrizations $\theta^{(0)}$ (including random initializations), the
algorithm consistently converged to the same $\theta^\star$ up to Monte
Carlo noise. This suggests that the pseudo-gradient/Pontryagin update
scheme is numerically stable in our test problems and finds at least a
locally optimal policy.

\section{Conclusion}
\label{sec:sec}
\subsection*{Summary of contributions}
We developed a framework for continuous-time stochastic control under partial information with discrete observations. From a theoretical perspective, we derived a Pontryagin maximum principle on the space of probability measures (belief space), yielding a coupled forward-backward system: a controlled filtering equation for the state's belief and an adjoint equation on the belief space. The adjoint's jump condition includes an extra term for the Bayes normalization, paralleling the correction term in nonlinear filtering, and under convexity conditions, we showed that the costate equals the (flat) gradient of the value functional along the optimal trajectory. Complementing this theory, our numerical methodology proposes a particle-based scheme that simulates the forward filter and approximates the adjoint using regression. By iteratively updating a parametric value/adjoint function via simulated data, the method finds near-optimal policies without solving high-dimensional PDEs.

Our results illustrate how an optimal POSOC controller balances information and cost. Between observation times, the control behaves like a full-information controller based on the current state estimate and at observation instants, it optimally trades off observation costs against the value of information. This approach leads to intuitive policies: for example, skipping or reducing observations when the state is well-estimated, or information is of low value, and increasing measurement frequency when uncertainty grows, or critical decisions are looming (thus actively managing the information-cost trade-off). The close agreement between our adjoint-based conditions and the dynamic programming principle, via the envelope relationship, gives confidence that these policies are indeed (locally) optimal.

\subsection*{Limitations} Our analysis assumes independent observation noise that does not directly enter the state dynamics or cost. Relaxing this, for instance allowing observation noise to influence the state, would require more advanced theory. The particle-regression algorithm yields only locally optimal solutions and can get stuck in suboptimal points if poorly initialized, which is a common issue in nonconvex optimization problems. Moreover, using a finite observation window or summary state (memory) is an approximation that may lose some long-term information, a commonplace trade-off between tractability and optimality. The method's computational cost can also grow with state dimension and horizon, due to the need for many particles and regression samples.

\subsection*{Future directions} In future research, we plan to co-design the filtering strategy with control. For example, optimizing the form of the observation summary (the mapping $\phi_n$ defining the memory state) alongside the policy. Extending the framework to continuous observation streams (recovering Kushner-Stratonovich filtering as a limit) or to event-triggered observations (where the controller decides when to observe) would broaden its applicability. Improving the efficiency and scalability of the particle-based solver is also crucial: using more expressive function approximators (e.g. deep neural networks) for the value function, and variance-reduction or parallelization techniques for the particle filter, could enable higher-dimensional applications. Rigorous convergence analysis of the algorithm and testing on complex real-world scenarios (e.g. multi-agent systems or adaptive experimental design) are other worthwhile steps. We believe that combining control and estimation in this way is a promising path toward tackling challenging partially observed decision problems. 
\newpage



\appendix

\appendix
\section{Proofs for Section~\ref{subsec:aux-optimality}}
\label{app:proofs-POSOC-aux}

\subsection{Proof of Proposition~\ref{prop:aux-opt}}
\label{app:proof-aux-opt}

We prove that any locally optimal policy in the partially observed class must satisfy the system \eqref{eq:aux_system}. The proof is by needle variations and is entirely pathwise.

\medskip\noindent
\textbf{Step 1: cost-to-go under a fixed policy.}
Fix $n\in\{0,\dots,N\}$, a realized observation history $\mf y^{[n]}=(y_{t_1},\dots,y_{t_n})$, and an admissible policy $u=(\alpha,\beta)\in\mcal U^{\mathrm{PO}}[t_n,T]$. For $t\in[t_n,t_{n+1})$ and $x\in\R^{d_x}$ define exactly as in the main text
\begin{equation}
	\label{eq:app-def-UR}
	U_{t}^{\mathrm R}(x,\mf y^{[n]}; u)
	:= \E\!\left[
	\int_t^T f_\tau(X_\tau,\alpha_\tau)\,d\tau
	+ \sum_{\{i:\,t_i\ge t\}} c_{t_i}(X_{t_i^-},\beta_{t_i})
	+ g(X_T)\ \Big|\ (X_t,\mf Y^{[n]})=(x,\mf y^{[n]})
	\right].
\end{equation}
On $(t_n,t_{n+1})$ the observation $\mf Y^{[n]}$ is frozen, so $(X_s)_{s\in[t_n,t_{n+1})}$ is a controlled diffusion with generator $\mcal G_{\alpha_s}$. By the Markov property,
\begin{equation}
	\label{eq:app-KB-UR}
	\frac{\partial}{\partial t} U_{t}^{\mathrm R}(x,\mf y^{[n]}; u)
	= - \mcal G_{\alpha_t} U_{t}^{\mathrm R}(x,\mf y^{[n]}; u)
	- f_t(x,\alpha_t),
	\qquad t\in[t_n,t_{n+1}),
\end{equation}
with terminal condition $U_T^{\mathrm R}(x,\mf y^{[N]}; u)=g(x)$.

At the observation time $t_{n+1}$, conditionally on $(X_{t_{n+1}^-},\mf y^{[n]})=(x,\mf y^{[n]})$, the new observation $Y_{t_{n+1}}$ has conditional density
\[
y \longmapsto \pi_{n+1}\bigl(y \mid x,\mf y^{[n]},\beta_{t_{n+1}}\bigr).
\]
Hence the jump relation for the pathwise cost is
\begin{equation}
	\label{eq:app-jump-UR}
	U_{t_{n+1}^-}^{\mathrm R}(x,\mf y^{[n]}; u)
	= c_{t_{n+1}}(x,\beta_{t_{n+1}})
	+ \int_{\R^{d_y}} U_{t_{n+1}}^{\mathrm R}\bigl(x,(\mf y^{[n]},y);\; u\bigr)\,
	\pi_{n+1}\bigl(y \mid x,\mf y^{[n]},\beta_{t_{n+1}}\bigr)\,dy.
\end{equation}
Equations \eqref{eq:app-KB-UR} and \eqref{eq:app-jump-UR} are the backward part of \eqref{eq:aux_system}, but stated for an \emph{arbitrary} policy $u$.

\medskip\noindent
\textbf{Step 2: a Itô-type identity.}
Let $u'=(\alpha',\beta')$ be another admissible policy. Consider the process
\[
s\longmapsto U_{s}^{\mathrm R}\bigl(X_s^{u'},\mf Y_{t_k}^{u'};\; u\bigr)
\quad\text{for } s\in[t_k,t_{k+1}),
\]
where $X^{u'}$ and $\mf Y^{u'}$ denote the state and the observations under $u'$. Applying Itô on every interval $[t_k,t_{k+1})$ and summing over $k$ yields the bookkeeping identity
\begin{align}
	J(u') - J(u)
	&= \sum_{k=0}^{N} \E\!\left[
	\int_{t_k}^{t_{k+1}}
	\Big(
	\mcal G_{\alpha'_s} U^{\mathrm R}_{s}(X^{u'}_s,\mf Y^{u'}_{t_k}; u)
	- \mcal G_{\alpha_s} U^{\mathrm R}_{s}(X^{u'}_s,\mf Y^{u'}_{t_k}; u)
	- f_s(X^{u'}_s,\alpha_s)
	+ f_s(X^{u'}_s,\alpha'_s)
	\Big)\,ds
	\right] \nonumber\\
	&\quad
	+ \sum_{k=1}^{N} \E\!\Big[
	U^{\mathrm R}_{t_k}(X^{u'}_{t_k},\mf Y^{u'}_{t_k}; u)
	- U^{\mathrm R}_{t_k^-}(X^{u'}_{t_k^-},\mf Y^{u'}_{t_k^-}; u)
	\Big]
	+ \sum_{k=1}^{N} \E\!\big[ c_{t_k}(X^{u'}_{t_k^-},\beta'_{t_k}) \big].
	\label{eq:app-Ito-UR}
\end{align}
This is the same identity as in the continuous-time verification lemma, except that we have an extra sum over the observation jumps.

\medskip\noindent
\textbf{Step 3: needle variation for the continuous control (optimality of $\alpha$).}
Let $\bar u=(\bar\alpha,\bar\beta)$ be a locally optimal policy and set
\[
\overline U_t(x,\mf y^{[n]})
:= U_t^{\mathrm R}(x,\mf y^{[n]}; \bar u).
\]
Fix $n$ and a history $\mf y^{[n]}$. Choose a small time interval
\[
I_\varepsilon = [\tilde t,\tilde t+\varepsilon]\subset (t_n,t_{n+1}),
\]
and pick any $\tilde\alpha\in\mcal A$.
Since on $(t_n,t_{n+1})$ the observation filtration is constant,
$\F_s^Y=\F_{t_n}^Y=\sigma(\mf Y^{[n]})$, the law of
$(X_s^{v},\mf Y^{[n]})$ always factors as
\[
\mcal{L}(X_s^{v},\mf Y^{[n]})
= \mu_s^{v}(dx \mid \mf y^{[n]}) \, \nu(d\mf y^{[n]}),
\qquad s\in(t_n,t_{n+1}),
\]
for any admissible $v$, where $\nu$ is the law of the observation history.

We now localize also in the observation variable. For $\varepsilon>0$ pick a measurable set
$B_\varepsilon^{(n)} \subset (\R^{d_y})^n$ such that $\nu(B_\varepsilon^{(n)})>0$ and
$B_\varepsilon^{(n)} \downarrow \{\mf y^{[n]}\}$ as $\varepsilon\downarrow0$. Define the
spike control
\[
\alpha_s^\varepsilon(\mf y^{[n]}) =
\begin{cases}
	\tilde\alpha, & (s,\mf y^{[n]}) \in I_\varepsilon \times B_\varepsilon^{(n)},\\
	\bar\alpha_s(\mf y^{[n]}), & \text{otherwise},
\end{cases}
\qquad
\beta_{t_k}^\varepsilon := \bar\beta_{t_k}\ \text{for all }k,
\]
and denote $u^\varepsilon:=(\alpha^\varepsilon,\bar\beta)$.

Plugging $(u^\varepsilon,\bar u)$ into \eqref{eq:app-Ito-UR}, and using that the two controls coincide outside $I_\varepsilon \times B_\varepsilon^{(n)}$, we obtain
\[
J(u^\varepsilon) - J(\bar u)
= \int_{I_\varepsilon} \int_{B_\varepsilon^{(n)}} \int_{\R^{d_x}}
\bigl(
f_s(x,\tilde\alpha)+\mcal G_{\tilde\alpha}\overline U_s(x,\mf y^{[n]})
- f_s(x,\bar\alpha_s)-\mcal G_{\bar\alpha_s}\overline U_s(x,\mf y^{[n]})
\bigr)
\,\mu_s^{u^\varepsilon}(dx \mid \mf y^{[n]})\,\nu(d\mf y^{[n]})\,ds.
\]
Divide both sides by $\varepsilon\,\nu(B_\varepsilon^{(n)})$ and let $\varepsilon\downarrow0$.
Since $I_\varepsilon \downarrow \{\tilde t\}$ and $B_\varepsilon^{(n)} \downarrow \{\mf y^{[n]}\}$,
and $s\mapsto \mu_s^{u^\varepsilon}(\cdot \mid \mf y^{[n]})$ is continuous on $(t_n,t_{n+1})$,
we get
\[
0 \le \int_{\R^{d_x}}
\bigl(
f_{\tilde t}(x,\tilde\alpha)+\mcal G_{\tilde\alpha}\overline U_{\tilde t}(x,\mf y^{[n]})
- f_{\tilde t}(x,\bar\alpha_{\tilde t})-\mcal G_{\bar\alpha_{\tilde t}}\overline U_{\tilde t}(x,\mf y^{[n]})
\bigr)\,
\overline\mu_{\tilde t}(dx \mid \mf y^{[n]}),
\]
where $\overline\mu_{\tilde t}(\cdot \mid \mf y^{[n]})$ is the belief induced by $\bar u$.
Since $\tilde\alpha\in\mcal A$ was arbitrary, this is equivalent to
\begin{equation}
	\label{eq:app-alpha-min}
	\bar\alpha_{\tilde t} \in \argmin_{\alpha\in\mcal A}
	\int_{\R^{d_x}}
	\bigl(f_{\tilde t}(x,\alpha)+\mcal G_{\alpha}\overline U_{\tilde t}(x,\mf y^{[n]})\bigr)
	\,\overline\mu_{\tilde t}(dx \mid \mf y^{[n]}),
\end{equation}

which is the third line of \eqref{eq:aux_system}.

\medskip\noindent
\textbf{Step 4: needle variation at an observation time (optimality of $\beta$).}
Fix $n$ and perturb only the control at the observation time $t_{n+1}$. Let $\tilde\beta\in\mcal B$ and define
\[
\beta_{t_{n+1}}^\varepsilon := \tilde\beta, \qquad
\beta_{t_k}^\varepsilon := \bar\beta_{t_k} \ (k\neq n+1), \qquad
\alpha^\varepsilon := \bar\alpha.
\]
Under $u^\varepsilon$, the dynamics coincide with those of $\bar u$ up to $t_{n+1}^-$, so
\[
\mcal{L}\bigl(X_{t_{n+1}^-},\mf Y^{[n]}\bigr) = \overline\mu_{t_{n+1}^-}(\cdot\mid \mf y^{[n]})\,\nu(d\mf y^{[n]}).
\]
Using \eqref{eq:app-jump-UR} with $u^\varepsilon$ and with $\bar u$ and plugging into \eqref{eq:app-Ito-UR}, the only nonzero contribution is at $t_{n+1}$:
\begin{align*}
	J(u^\varepsilon) - J(\bar u)
	&= \E\bigl[c_{t_{n+1}}(X_{t_{n+1}^-},\tilde\beta) - c_{t_{n+1}}(X_{t_{n+1}^-},\bar\beta_{t_{n+1}})\bigr]\\
	&\quad + \E\biggl[
	\int_{\R^{d_y}}
	\overline U_{t_{n+1}}\bigl(X_{t_{n+1}^-},(\mf Y^{[n]},y)\bigr)
	\Big(
	\pi_{n+1}(y\mid X_{t_{n+1}^-},\mf Y^{[n]},\tilde\beta)
	- \pi_{n+1}(y\mid X_{t_{n+1}^-},\mf Y^{[n]},\bar\beta_{t_{n+1}})
	\Big)\,dy
	\biggr].
\end{align*}
Condition on $(X_{t_{n+1}^-},\mf Y^{[n]})=(x,\mf y^{[n]})$ and integrate w.r.t.\ $\overline\mu_{t_{n+1}^-}(dx\mid \mf y^{[n]})\nu(d\mf y^{[n]})$ to obtain
\begin{align*}
	J(u^\varepsilon) - J(\bar u)
	&= \int_{(\R^{d_y})^n} \int_{\R^{d_x}}
	\Bigg(
	c_{t_{n+1}}(x,\tilde\beta)
	+ \int_{\R^{d_y}} \overline U_{t_{n+1}}(x,(\mf y^{[n]},y))\,
	\pi_{n+1}(y\mid x,\mf y^{[n]},\tilde\beta)\,dy
	\Bigg) \overline\mu_{t_{n+1}^-}(dx\mid \mf y^{[n]})\,\nu(d\mf y^{[n]})\\
	&\quad-
	\int_{(\R^{d_y})^n} \int_{\R^{d_x}}
	\Bigg(
	c_{t_{n+1}}(x,\bar\beta_{t_{n+1}})
	+ \int_{\R^{d_y}} \overline U_{t_{n+1}}(x,(\mf y^{[n]},y))\,
	\pi_{n+1}(y\mid x,\mf y^{[n]},\bar\beta_{t_{n+1}})\,dy
	\Bigg) \overline\mu_{t_{n+1}^-}(dx\mid \mf y^{[n]})\,\nu(d\mf y^{[n]}).
\end{align*}
Local optimality of $\bar u$ implies this difference is $\ge 0$ for every $\tilde\beta\in\mcal B$, hence for every fixed history $\mf y^{[n]}$ and for $\overline\mu_{t_{n+1}^-}(\cdot\mid \mf y^{[n]})$–a.e.\ $x$,
\[
\bar\beta_{t_{n+1}}
\in \argmin_{\beta\in\mcal B}
\Bigg\{
c_{t_{n+1}}(x,\beta)
+ \int_{\R^{d_y}} \overline U_{t_{n+1}}(x,(\mf y^{[n]},y))\,
\pi_{n+1}(y\mid x,\mf y^{[n]},\beta)\,dy
\Bigg\},
\]
and averaging w.r.t.\ $\overline\mu_{t_{n+1}^-}(\cdot\mid \mf y^{[n]})$ gives precisely the last line of \eqref{eq:aux_system}.

\medskip\noindent
\textbf{Step 5: forward evolution of the belief.}
Between observation dates the conditional law of $X_t$ given $\mf y^{[n]}$ under $\bar u$ satisfies the Fokker--Planck equation
\[
\dot{\overline\mu}_t(\cdot\mid \mf y^{[n]})
= \mcal G_{\bar\alpha_t}^{*}\,\overline\mu_t(\cdot\mid \mf y^{[n]}),
\qquad t\in(t_n,t_{n+1}),
\]
and at $t_{n+1}$ the Bayesian update with the density $\pi_{n+1}$ gives
\[
\overline{\mu}_{t_{n+1}}(x\mid \mf y^{[n]},y)
= \frac{\pi_{n+1}(y\mid x,\mf y^{[n]},\bar\beta_{t_{n+1}})}
{L_{n+1}(y;\mf y^{[n]},\bar\beta_{t_{n+1}})}\,
\overline{\mu}_{t_{n+1}^-}(x\mid \mf y^{[n]}),
\]
which are exactly the fourth and fifth lines of \eqref{eq:aux_system}.

\medskip
Collecting \eqref{eq:app-KB-UR}, \eqref{eq:app-jump-UR}, the minimization conditions \eqref{eq:app-alpha-min} and the discrete-time minimization above, together with the forward belief dynamics, we obtain the full auxiliary  system \eqref{eq:aux_system}. This completes the proof.
\qedhere
\subsection{Proof of Proposition~\ref{prop:FB-belief-cases}}
\label{app:proof-FB-belief}

\begin{proof}
	We prove that saddle points of
	\[
	\inf_{u,\mu}\sup_{\lambda}\;\mcal L(u,\mu,\lambda)
	\]
	satisfy the system \eqref{eq:FOC-belief}. The key points are:
	(i) the cost is linear in the belief,
	(ii) between observation times the observation filtration is constant,
	(iii) at observation times the constraint is a Bayesian update, so the variation must be taken conditionally.
	
	\medskip
	\noindent\textbf{1. Linear-in-belief form of the cost.}
	By assumption,
	\[
	J(u)
	= \E\!\bigg[\int_0^T \!\!\lrang{f_t(\cdot,\alpha_t)}{\mu_t}\dl t
	+ \sum_{n=1}^N \lrang{c_{t_n}(\cdot,\beta_n)}{\mu_{t_n^-}}
	+ \lrang{g(\cdot)}{\mu_T}\bigg],
	\]
	where, on each interval \((t_n,t_{n+1})\) we write
	\(\mu_t(\cdot) \equiv \mu_t(\cdot\mid \mf Y^{[n]})\).
	The belief \(\mu\) is constrained by the controlled filtering dynamics
	\[
	\dot\mu_t = \mcal G_{\alpha_t}^* \mu_t, \quad t\in(t_n,t_{n+1}), 
	\qquad
	\mu_{t_{n+1}} = \mcal K_{\beta_{n+1},Y_{t_{n+1}}}\!\big(\cdot;\,\mu_{t_{n+1}^-}\big).
	\]
	
	\medskip
	\noindent\textbf{2. Lagrangian.}
	Introduce an adjoint (costate) \(\lambda_t(\cdot,\mf Y^{[n]})\) which, on each interval \((t_n,t_{n+1})\), is \(\F_{t_n}^Y\)-measurable.
	The Lagrangian is
	\begin{align*}
		\mcal L(u,\mu,\lambda)
		&= J(u)
		- \E\!\left[\sum_{n=1}^{N+1}\int_{t_{n-1}}^{t_n}
		\left\langle \lambda_t^{[n]}(\cdot,\mf Y^{[n]}),\,
		\dot\mu_t^{[n]}(\cdot\mid \mf Y^{[n]})-\mcal G_{\alpha_t}^*\,\mu_t^{[n]}(\cdot\mid \mf Y^{[n]})
		\right\rangle \dl t\right]\\
		&\quad - \E\!\left[\sum_{n=1}^{N}
		\left\langle \lambda_{t_n}^{[n]}(\cdot,\mf Y^{[n]}),\,
		\mu_{t_n}^{[n]}(\cdot\mid \mf Y^{[n]})-\mcal K_{\beta_{t_n},Y_{t_n}}\!\big(\cdot;\,\mu_{t_n^-}^{[n-1]}(\cdot\mid \mf Y^{[n-1]})\big)
		\right\rangle \right].
	\end{align*}
	(Here \(\lrang{\cdot}{\cdot}\) is the duality between bounded test functions and finite measures on \(\R^{d_x}\).)
	
	A saddle point \((\bar u,\bar\mu,\bar\lambda)\) must satisfy that the first variation in each direction vanishes.
	
	\medskip
	\noindent\textbf{3. Variation w.r.t. \(\lambda\): recovery of the belief flow.}
	Let \(\delta\lambda\) be an admissible variation, i.e. on \((t_n,t_{n+1})\) it is \(\F_{t_n}^Y\)-measurable. Since \(\mcal L\) is affine in \(\lambda\),
	\[
	\delta_\lambda \mcal L
	= - \E\!\left[\sum_{n=1}^{N+1}\int_{t_{n-1}}^{t_n}
	\big\langle \delta\lambda_t^{[n]},\,\dot\mu_t^{[n]}-\mcal G_{\alpha_t}^*\mu_t^{[n]}\big\rangle \dl t\right]
	- \E\!\left[\sum_{n=1}^{N}
	\big\langle \delta\lambda_{t_n}^{[n]},\,\mu_{t_n}^{[n]}-\mcal K_{\beta_{t_n},Y_{t_n}}(\cdot;\mu_{t_n^-}^{[n-1]})\big\rangle \right].
	\]
	Localizing in time and in the observation history as in the previous subsection, we can choose
	\(\delta\lambda\) supported in an arbitrarily small time/observation tube;
	hence each bracket must vanish \emph{pathwise}, i.e. for \(\nu\)-a.e. history \(\mf y^{[n]}\):
	\[
	\dot\mu_t(\cdot\mid \mf y^{[n]}) = \mcal G_{\alpha_t}^* \mu_t(\cdot\mid \mf y^{[n]}),
	\quad t\in(t_n,t_{n+1}),
	\]
	and
	\[
	\mu_{t_{n+1}}(\cdot\mid \mf y^{[n]},y)
	= \mcal K_{\beta_{t_{n+1}},y}\!\big(\cdot;\,\mu_{t_{n+1}^-}(\cdot\mid \mf y^{[n]})\big),
	\quad \text{for }\pi_{n+1}(\cdot\mid x,\mf y^{[n]},\beta_{t_{n+1}})\text{-a.e. }y.
	\]
	This recovers the forward (filtering) part of \eqref{eq:FOC-belief}.
	
	\medskip
	\noindent\textbf{4. Variation w.r.t. \(\mu\) on \((t_n,t_{n+1})\): backward adjoint.}
	Now fix \(n\) and a history \(\mf y^{[n]}\). On \((t_n,t_{n+1})\) the Lagrangian contains, after expanding \(J(u)\),
	\[
	\int_{t_n}^{t_{n+1}} \Big(
	\lrang{f_t(\cdot,\alpha_t)}{\mu_t}
	- \lrang{\lambda_t}{\dot\mu_t - \mcal G_{\alpha_t}^* \mu_t}
	\Big)\dl t.
	\]
	Integrating by parts the term \(\lrang{\lambda_t}{\dot\mu_t}\) on \((t_n,t_{n+1})\) we get
	\[
	\int_{t_n}^{t_{n+1}} \!\!\lrang{f_t(\cdot,\alpha_t) + \mcal G_{\alpha_t}\lambda_t}{\mu_t}\,\dl t
	+ \lrang{\lambda_{t_n^-}}{\mu_{t_n^-}} - \lrang{\lambda_{t_{n+1}^-}}{\mu_{t_{n+1}^-}}.
	\]
	(We used that \(\lrang{\lambda_t}{\mcal G_{\alpha_t}^*\mu_t} = \lrang{\mcal G_{\alpha_t}\lambda_t}{\mu_t}\).)
	Perturb now \(\mu\) by an arbitrary predictable \(\delta\mu_t(\cdot\mid \mf y^{[n]})\) on \((t_n,t_{n+1})\). The variation is
	\[
	\int_{t_n}^{t_{n+1}} \!\!\lrang{
		f_t(\cdot,\alpha_t) + \mcal G_{\alpha_t}\lambda_t + \dot\lambda_t
	}{\delta\mu_t}\,\dl t.
	\]
	By localization (we can choose \(\delta\mu_t\) supported in a small time interval and at a fixed history),
	this implies, for a.e. \(t\in(t_n,t_{n+1})\) and every fixed \(\mf y^{[n]}\),
	\[
	\dot\lambda_t(x,\mf y^{[n]})
	= -\big(\mcal G_{\alpha_t}\lambda_t(x,\mf y^{[n]}) + f_t(x,\alpha_t)\big),
	\]
	which is exactly the second line of \eqref{eq:FOC-belief} (with \(\bar\alpha_t\) when evaluated at the optimal control).
	
	\medskip
	\noindent\textbf{5. Variation w.r.t. \(\mu\) at an observation time (jump of \(\lambda\)).}
	The only terms involving \(\mu_{t_{n+1}^-}(\cdot\mid \mf y^{[n]})\) in the Lagrangian are
	\[
	\lrang{c_{t_{n+1}}(\cdot,\beta_{t_{n+1}})}{\mu_{t_{n+1}^-}}
	+ \E\big[\,\lrang{\lambda_{t_{n+1}}}{\mcal K_{\beta_{t_{n+1}},Y_{t_{n+1}}}(\cdot;\mu_{t_{n+1}^-})}\,\big|\F_{t_n}^Y\big]
	- \lrang{\lambda_{t_{n+1}^-}}{\mu_{t_{n+1}^-}}.
	\]
	Work pathwise: fix \(\mf y^{[n]}\) and denote by
	\[
	L_{n+1}(y;\mf y^{[n]},\beta)
	:= \int_{\R^{d_x}} \pi_{n+1}(y\mid x',\mf y^{[n]},\beta)\,
	\mu_{t_{n+1}^-}(\dl x'\mid \mf y^{[n]})
	\]
	the predictive density of \(Y_{t_{n+1}}\) under \(\mu_{t_{n+1}^-}\). Then
	\[
	\mcal K_{\beta,y}(\cdot;\mu_{t_{n+1}^-})
	= \frac{\pi_{n+1}(y\mid \cdot,\mf y^{[n]},\beta)}{L_{n+1}(y;\mf y^{[n]},\beta)}\,
	\mu_{t_{n+1}^-}(\cdot\mid \mf y^{[n]}).
	\]
	Hence
	\begin{align*}
		&\E\big[\,\lrang{\lambda_{t_{n+1}}}{\mcal K_{\beta_{t_{n+1}},Y_{t_{n+1}}}(\cdot;\mu_{t_{n+1}^-})}\,\big|\F_{t_n}^Y\big]\\
		&= \int_{\R^{d_y}}\!\int_{\R^{d_x}}\!
		\lambda_{t_{n+1}}(x,(\mf y^{[n]},y))\,
		\frac{\pi_{n+1}(y\mid x,\mf y^{[n]},\beta_{t_{n+1}})}{L_{n+1}(y;\mf y^{[n]},\beta_{t_{n+1}})}\,
		\mu_{t_{n+1}^-}(\dl x\mid \mf y^{[n]})\,L_{n+1}(y;\mf y^{[n]},\beta_{t_{n+1}})\,\dl y\\
		&= \int_{\R^{d_x}}\!\mu_{t_{n+1}^-}(\dl x\mid \mf y^{[n]})
		\int_{\R^{d_y}}\!\lambda_{t_{n+1}}(x,(\mf y^{[n]},y))
		\pi_{n+1}(y\mid x,\mf y^{[n]},\beta_{t_{n+1}})\,\dl y.
	\end{align*}
	Now let \(\delta\mu_{t_{n+1}^-}\) be any \(\F_{t_n}^Y\)-measurable variation. The first-order condition in \(\mu_{t_{n+1}^-}\) reads
	\begin{align*}
		0
		&= \int_{\R^{d_x}} \Bigg[
		c_{t_{n+1}}(x,\beta_{t_{n+1}})
		+ \int_{\R^{d_y}}\!\lambda_{t_{n+1}}(x,(\mf y^{[n]},y))\,
		\pi_{n+1}(y\mid x,\mf y^{[n]},\beta_{t_{n+1}})\,\dl y
		- \lambda_{t_{n+1}^-}(x,\mf y^{[n]}) \\
		&\qquad\quad
		- \int_{\R^{d_x}}\!\mu_{t_{n+1}^-}(\dl x'\mid \mf y^{[n]})
		\int_{\R^{d_y}}\!\lambda_{t_{n+1}}(x',(\mf y^{[n]},y))\,
		\frac{\pi_{n+1}(y\mid x',\mf y^{[n]},\beta_{t_{n+1}})\,\pi_{n+1}(y\mid x,\mf y^{[n]},\beta_{t_{n+1}})}{L_{n+1}(y;\mf y^{[n]},\beta_{t_{n+1}})}\,\dl y
		\Bigg]\,
		\delta\mu_{t_{n+1}^-}(\dl x).
	\end{align*}
	Since \(\delta\mu_{t_{n+1}^-}\) is arbitrary, the bracket must vanish, and we obtain exactly the jump relation in \eqref{eq:FOC-belief}:
	\begin{align*}
		\lambda_{t_{n+1}^-}(x,\mf y^{[n]})
		&= c_{t_{n+1}}\!\left(x,\beta_{t_{n+1}}\right)
		+ \int_{\R^{d_y}}\!\lambda_{t_{n+1}}(x,(\mf y^{[n]},y))\,\pi_{n+1}(y\mid x,\mf y^{[n]},\beta_{t_{n+1}})\,\dl y \\
		&\quad - \int_{\R^{d_x}}\!\mu_{t_{n+1}^-}(\dl x'\mid \mf y^{[n]})
		\int_{\R^{d_y}}\!\lambda_{t_{n+1}}(x',(\mf y^{[n]},y))\,
		\frac{\pi_{n+1}(y\mid x',\mf y^{[n]},\beta_{t_{n+1}})\,\pi_{n+1}(y\mid x,\mf y^{[n]},\beta_{t_{n+1}})}{L_{n+1}(y;\mf y^{[n]},\beta_{t_{n+1}})}\,\dl y.
	\end{align*}
	The terminal condition \(\lambda_T(x,\mf y^{[N]})=g(x)\) comes from the variation of the terminal term \(\lrang{g}{\mu_T}\).
	
	\medskip
	\noindent\textbf{6. Variation w.r.t. \(\alpha\): continuous control optimality.}
	On \((t_n,t_{n+1})\) the \(\alpha\)-dependent part of \(\mcal L\) is
	\[
	\int_{t_n}^{t_{n+1}} \lrang{f_t(\cdot,\alpha_t) + \mcal G_{\alpha_t}\lambda_t}{\mu_t}\,\dl t.
	\]
	Let \(\delta\alpha_t\) be any \(\F_t^Y\)-measurable perturbation, supported in a small subinterval.
	Then
	\[
	0 = \delta_\alpha \mcal L
	= \int_{t_n}^{t_{n+1}} \int_{\R^{d_x}}
	\Big(\partial_\alpha f_t(x,\bar\alpha_t)
	+ \partial_\alpha \mcal G_{\alpha}\lambda_t(x,\mf y^{[n]})\big|_{\alpha=\bar\alpha_t}\Big)
	\,\bar\mu_t(\dl x\mid \mf y^{[n]})\,\delta\alpha_t \,\dl t.
	\]
	Since we can localize in \(t\) and \(\mf y^{[n]}\) and take arbitrary \(\delta\alpha_t\), this is equivalent to the variational inequality
	\[
	\bar\alpha_t \in \argmin_{\alpha\in\mcal A}
	\int_{\R^{d_x}} \big(f_t(x,\alpha)+\mcal G_{\alpha}\lambda_t(x,\mf y^{[n]})\big)\,
	\bar\mu_t(\dl x\mid \mf y^{[n]}),
	\quad t\in(t_n,t_{n+1}),
	\]
	which is the third line of \eqref{eq:FOC-belief}.
	
	\medskip
	\noindent\textbf{7. Variation w.r.t. \(\beta\): discrete control optimality.}
	At time \(t_{n+1}\) the \(\beta\)-dependent part of \(\mcal L\) (conditioned on \(\mf y^{[n]}\)) is
	\[
	\int_{\R^{d_x}} c_{t_{n+1}}(x,\beta)\,\mu_{t_{n+1}^-}(\dl x\mid \mf y^{[n]})
	+ \int_{\R^{d_x}} \mu_{t_{n+1}^-}(\dl x\mid \mf y^{[n]})
	\int_{\R^{d_y}} \lambda_{t_{n+1}}(x,(\mf y^{[n]},y))\,\pi_{n+1}(y\mid x,\mf y^{[n]},\beta)\,\dl y.
	\]
	Let \(\delta\beta\) be \(\F_{t_{n+1}^-}^Y\)-measurable. Localizing in \(\mf y^{[n]}\) we obtain the pointwise condition
	\[
	\bar\beta_{n+1}\in\argmin_{\beta\in\mcal B}\!
	\int_{\R^{d_x}}\!\Bigg(c_{t_{n+1}}(x,\beta)
	+ \int_{\R^{d_y}}\!\lambda_{t_{n+1}}(x,(\mf y^{[n]},y))\,\pi_{n+1}(y\mid x,\mf y^{[n]},\beta)\,\dl y\Bigg)
	\,\mu_{t_{n+1}^-}(\dl x\mid \mf y^{[n]}),
	\]
	which is the last line of \eqref{eq:FOC-belief}.
	
	\medskip
	Putting together Steps 3–7 yields exactly the KKT/PMP system \eqref{eq:FOC-belief}.
\end{proof}


\footnotesize
\bibliographystyle{acm} 
\bibliography{references}    

\end{document}